\documentclass[11pt,a4paper]{amsart} 
\usepackage[english]{babel}
\usepackage{amssymb,xspace,enumerate,graphics}
\usepackage{amstext}
\theoremstyle{plain}
\usepackage{amsbsy,amssymb,amsfonts,latexsym}

\date{\today}
\title{Extension of holomorphic functions  defined on
singular analytic spaces with growth estimates}
\author{William ALEXANDRE}
\address{Laboratoire Paul Painlev\'e U.M.R. CNRS 8524, U.F.R. de
Math\'ematiques,  cit\'e scientifique, Universit\'e Lille 1, F59 655 Villeneuve d'Ascq Cedex, France.}
\email{ william.alexandre@math.univ-lille1.fr}
\thanks{The first author is partially supported by A.N.R. BL-INTER09-CRARTIN}
\author{Emmanuel MAZZILLI}
\address{Laboratoire Paul Painlev\'e U.M.R. CNRS 8524, U.F.R. de
Math\'ematiques,  cit\'e scientifique, Universit\'e Lille 1, F59 655 Villeneuve d'Ascq Cedex, France.}
\email{ emmanuel.mazzilli@math.univ-lille1.fr}
\subjclass[2000]{32A22, 32A26, 32A27, 32A37, 32A40, 32A55,\\32C30, 32D15}
\keywords{Analytic spaces, holomorphic extensions, residue currents, integral representations.}
\date{}
\newtheorem{theorem}{Theorem}[section]
\newtheorem{lemma}[theorem]{Lemma}
\newtheorem{fact}[theorem]{Fact}
\newtheorem{example}[theorem]{Example}
\newtheorem{proposition}[theorem]{Proposition}
\newtheorem{corollary}[theorem]{Corollary}
\newtheorem{remark}{\it Remark}
\newtheorem{definition}[theorem]{Definition}
\def \pint {\vbox{ \hbox to 5 pt {\hfil \vrule height 4pt}\hrule}\hskip 3pt}
\def\leqs{\lesssim}
\def\geqs{\gtrsim}
\def\eqs{\eqsim}
\def\cc{\mathbb{C}}
\def\rr{\mathbb{R}}
\def\nn{\mathbb{N}}
\def\zz{\mathbb{Z}}

\def \qed {\hbox{\hskip 5pt} \vbox{\hrule \hbox to 5pt 
{\vrule height 4.2pt \hfil \vrule}\hrule}}
\def \pint {\vbox{ \hbox to 5 pt {\hfil \vrule height 4pt}\hrule}\hskip 3pt}

\newcommand{\cal}{\mathcal}
\newcommand{\wrt}{with respect to }

\newcommand{\app}[5]
{#1: \left \{ 
\begin {array}{ccl} 
#2 & \longrightarrow & #3 \\#4 & \longmapsto & #5 \end {array} \right .}

\newcommand{\tauzv}{\tau\left(z,v,|\rho(z)|\right)}
\newcommand{\diffp}[2]{\frac{\partial #1}{\partial #2}}
\newcommand{\mlabel}[1]{\label {#1}}
\renewcommand{\over}[2]{\genfrac{}{}{0pt}{}{#1}{#2}}

\newcommand{\p}[2]{{\cal P}_{#1}(#2)}

\newcommand{\pk}[1]{{\cal P}_{\kappa|\rho(z_{#1})|}(z_{#1})}

\newcommand{\pr} {\noindent{\it Proof:} }
\newcommand{\ko} {Koranyi }
\newcommand{\kb}[2] {{\cal P}_{#1}(#2)}

\newcommand{\oo}{{\cal O}}
\begin{document}
\pagestyle{plain}

\begin{abstract}
Let $D$ be a strictly convex domain and $X$ be an analytic subset of $\cc^n$ such that $X\cap D\neq \emptyset$ and $X\cap bD$ is transverse. We first  give necessary conditions for a function holomorphic on $D\cap X$ to admit  a holomorphic extension belonging to $L^q(D),$ $q\in [1,+\infty]$. When $n=2$ and $q<+\infty$, we then prove that this condition is also sufficient. When $q=+\infty$ we prove that this condition implies the existence of a $BMO$-holomorphic extension. In both cases, the extensions are given by mean of integral representation formulas and new residue currents.\end{abstract}
\maketitle
\section{Introduction}
In the last few years, many researches have been done on classical problems in complex analysis in the case of singular spaces; for example the $\overline\partial$-Neumann operator has been studied in \cite{Rup2} by Ruppenthal, the Cauchy-Riemann equation in \cite{AS,DFV,ForII,Rup0,Rup1} by Andersson, Samuelsson, Diederich, Forn\ae ss, Vassiliadou, Ruppenthal, ideals of holomorphic functions on analytic spaces in \cite{ASS} by Andersson, Samuelsson and Sznajdman,  problems of extensions and restrictions of holomorphic functions on analytic spaces in \cite{DiMa0, Duquenoy} by Diederich, Mazzilli and Duquenoy.

In this article, we will be interested in problems of extension of holomorphic functions defined on an analytic space. Let $D$ be a bounded pseudoconvex domain of $\cc^n$ with smooth boundary, let $f$ be a holomorphic function in a neighborhood of $D$ and let $X=\{z, f(z)=0\}$ be an analytic set such that $D\cap X\neq \emptyset$. The first extension problem that one can consider is the following one~: Is it true that a function $g$ which is holomorphic on $D\cap X$ has a holomorphic extension on $D$ ?\\
It is known by Cartan's theorem B that the answer to this question is affirmative and that any function $g$ holomorphic on $X\cap D$ has a holomorphic extension $G$ on the whole domain $D$ if and only if $D$ is pseudoconvex. More difficulties arise when we ask $G$ to satisfy some growth conditions like being in $L^q(D)$ or in $BMO(D)$. This question has been widely studied by many authors under different assumptions on $D$ or $X$. In \cite{Ohsawa}, Ohsawa-Takegoshi proved when $X$ is a hyperplane that any $g\in L^2(X\cap D)\cap {\cal O}(X\cap D)$ admits an extension $G\in L^2(D)\cap \oo (D)$. This result was generalized to the case of manifolds of higher codimension in \cite{OhsawaII} by Ohsawa. In \cite{Berndtsson}, Berndtsson investigated the case of singular varieties and got a condition on $g$ which implies that it admits a holomorphic $L^2$ extension on $D$. However this condition requires that $g$ vanishes on the singularities of $X$ and thus $g\equiv 1$ does not satisfy this condition while it can trivially be extended holomorphically.\\
Assuming that $D$ is strictly pseudoconvex and that $X$ is a manifold, Henkin proved in \cite{Henkin} that any $g\in L^\infty(D\cap X)\cap\
\oo(D\cap X)$ has an extension in $L^\infty(D)\cap\oo(D)$, provided that $bD$, the boundary of $D$, and $X$ are in general position. Cumenge in \cite{Cum} generalized this result to the case of Hardy spaces and Amar in \cite{Ama} removed the hypothesis of general position of $bD$ and $X$ assumed in \cite{Henkin}. The case of $L^\infty$ extensions has also been investigated in the case of weak (pseudo)convexity. In \cite{DiMa2} Diederich and Mazzilli proved that when $D$ is convex of finite type and $X$ is a hyperplane, any $g\in L^\infty(D\cap X)\cap\oo(D\cap X)$ is the restriction of some $G\in L^\infty(D)\cap\oo(D)$. In \cite{WA}, again for $D$ convex of finite type but for $X$ a manifold, a sufficient and nearly necessary condition on $X$ was given under which any function $g$ which is bounded and holomorphic on $X\cap D$ is the restriction of a bounded holomorphic function on $D$. This restriction problem was also studied in \cite{Jasiczak} by Jasiczak for $D$ a pseudoconvex domain of finite type in $\cc^2$ and $X$ a manifold.

In this article we consider a strictly convex domain $D$ of $\cc^n$ and an analytic subset $X$ of $\cc^n$ such that $X\cap D\neq\emptyset$ and $X\cap bD$ is transverse in the sense of tangent cones. We give necessary conditions and, when $n=2$,  sufficient conditions under which a function $g$ holomorphic in $X\cap D$ admits a holomorphic extension in the class $BMO(D)$ or $L^q(D)$, $q\in [1,+\infty)$.
\par\medskip
Let us write $D$ as $D=\{z\in\cc^n,\ \rho(z)<0\}$ where $\rho$ is a smooth strictly convex function defined on $\cc^n$ such that the gradient of $\rho$ does not vanish in a neighborhood $\cal U$ of $bD$. We denote by $D_r$, $r\in\rr$, the set $D_r=\{z\in\cc^n,\ \rho(z)<0\}$, by $\eta_\zeta$  the outer unit normal to $bD_{\rho(\zeta)}$ at a point $\zeta\in{\cal U}$ and by $v_\zeta$ a smooth complex tangent vector field at $\zeta$ to $bD_{\rho(\zeta)}$. Our first result is the following.
\begin{theorem}\label{th0}
For $n=2$, there exists two integers $k,l\geq 1$ depending only from $X$ such that if $g$ is a holomorphic function on $X\cap D$ which has a $C^\infty$ smooth extension $\tilde g$ on $D$ which satisfies
\begin{enumerate}[(i)]
 \item \label{th0i} there exists $N\in\nn$ such that $|\rho|^N \tilde g$ vanishes to order $l$ on $bD$,
 \item \label{th0ii}there exists $q\in[1,+\infty]$ such that  $\left|\diffp{^{\alpha+\beta}\tilde g}{\overline{\eta_\zeta}^\alpha\partial \overline{v_\zeta}^\beta}\right||\rho|^{\alpha+\frac\beta 2}$ belongs to $L^q(D)$ for all non-negative integers $\alpha$ and $\beta$  with $\alpha+\beta\leq k$,
 \item\label{th0iii} $\diffp{^{\alpha+\beta}\tilde g}{\overline{\eta_\zeta}^\alpha\partial \overline{v_\zeta}^\beta}=0$ on $X\cap D$ for all non-negative integers $\alpha$ and $\beta$  with $\alpha+\beta\leq k$,
\end{enumerate}
then $g$ has a holomorphic extension $G$ in $L^q(D)$ when $q<+\infty$ and in $BMO(D)$ when $q=+\infty$. Moreover, up to a uniform multiplicative constant  depending only from $k$, $l$ and $N$, the norm of $G$ is  bounded 
by the supremum of the $L^q$-norm of $\zeta\mapsto \left|\diffp{^{\alpha+\beta}\tilde g}{\overline{\eta_\zeta}^\alpha\partial \overline{v_\zeta}^\beta}(\zeta)\right||\rho(\zeta)|^{\alpha+\frac\beta 2}$ for $\alpha ,\beta $ with $\alpha+\beta\leq k$.
\end{theorem}
In  Lemma \ref{lemma2}, Corollary \ref{th2} and Theorem \ref{th5}, we will give conditions under which a function $g$ holomorphic on $X\cap D$ admits a smooth extension on $D$ which satisfies the assumption of Theorem \ref{th0}.\\
Let us mention that the integer $k$ in Theorem \ref{th0} is in fact equal to the maximum of the order of the singularities of $X$ and the hypothesis of Theorem \ref{th0} can be relaxed a little in the following way. The theorem is still valid if for all singularities $z_0\in X\cap\overline D$ of $X$ of order $k_0$, we check the hypothesis (\ref{th0ii}) and (\ref{th0iii}) with $k$ replaced by $k_0$ and $D$ replaced by ${\cal U}_0\cap D$ where ${\cal U}_0$ is a neighborhood of $z_0$.
\par\smallskip
The holomorphic extension of Theorem \ref{th0} is given by an integral operator combining the Berndtsson-Andersson reproducing kernel and a residue current. In \cite{Ama}, Amar pointed out for the first time the importance of the current $\overline\partial\left[ \frac1f\right]$ in the problem of extension. In \cite{Duquenoy} the extension is given by an operator constructed by Passare which uses the classical residue current $\overline\partial \left[ \frac1f\right]$  (see \cite{Passare}). However, as pointed out in \cite{Duquenoy}, it is not so easy to handle the case of singularities of order greater than 2 and the classical currents do not give a good extension in this case. To overcome this difficulty we  have to adapt a construction due to the second author of new residue currents which will play the role of $\overline\partial\left[ \frac1f\right]$ (see  \cite{Maz1} and \cite{Mazzilli2}). The extension given by Theorem \ref{th0} will be obtained via a linear operator which uses a Berndtsson-Andersson reproducing kernel and these new currents (see Section \ref{section3}).
\par\smallskip
Observe that in Theorem \ref{th0} we assume the existence of a smooth extension $\tilde g$ satisfying properties (\ref{th0i}), (\ref{th0ii}) and (\ref{th0iii}), whereas no such assumption is made in the previous articles we quoted, which deal with extension problems. It should be pointed out that while boundedness is a sufficient hypothesis in order to obtain a bounded holomorphic extension when $X$ is a manifold (see \cite{ WA,Ama,Cum,DiMa2}), it is not possible to obtain $L^\infty$ or even $L^2$ extensions when $X$ has singularities if we only assume that $g$ is bounded on $X\cap D$ (see \cite{DiMa0}) : a stronger condition is needed. Actually, even if in the manifold case no smooth extension is assumed to exist, a smooth extension, which satisfies (\ref{th0ii}) and (\ref{th0iii}), is constructed for example in \cite{Cum,DiMa2,WA}. This is done as follows. When $X$ is a manifold, let us locally write $X$ as $X=\{(z',\alpha(z')),\ z'\in\cc^{n-1}\}$, with $\alpha$ holomorphic. If for $z=(z_1,\ldots, z_n)$ we set $z'=(z_1,\ldots, z_{n-1})$, then the function $\tilde g$ defined by  $\tilde g(z):=g(z',\alpha(z'))$ is a local holomorphic extension of $g$. Gluing all these local extensions together we get a smooth extension which will satisfy (\ref{th0ii}) and (\ref{th0iii}).  In some sense, the way the local holomorphic extension is constructed in the manifold case is a kind of interpolation : $\tilde g(z',\cdot)$ is the polynomial of degree $0$ which interpolates $g(z',\alpha(z'))$ at the point $z_n=\alpha(z')$. 
Following this idea, we will construct in Section \ref{section5} a local holomorphic extension by interpolation. Provided we have a good control of the polynomials which interpolate $g$ on the different sheets of $X$, gluing together these local extensions, we will obtain an appropriate smooth extension. The control of the interpolating polynomials will be achieved thanks to an assumption on the divided differences we can build with $g$ between the different sheets of $X$. This will give us simple numerical conditions under which the function $g$ has a smooth extension $\tilde g$ which satisfies (\ref{th0i}), (\ref{th0ii}) and (\ref{th0iii}) from Theorem \ref{th0} (see Theorem \ref{th2} and \ref{th5}).
The divided differences are defined has follows.
\par\smallskip 
For $z\in D$, $v$ a unit vector in $\cc^n$, and $\varepsilon$ a positive real number we set $\Delta_{z,v}(\varepsilon)=\{z+\lambda v,\ |\lambda|<\varepsilon\}$ and 
$$\tau(z,v,\varepsilon)=\sup\{\tau>0,\ \rho(z+\lambda v)-\rho(z)<\varepsilon\text{ for all } \lambda\in\cc,\ |\lambda|<\tau\}.$$
Therefore $\tau(z,v,\varepsilon)$ is the maximal radius $r>0$ such that the disc $\Delta_{z,v}\left(r\right)$  is in $D_{\rho(z)+\varepsilon}$. It is also the distance from $z$ to $bD_{\rho(z)+\varepsilon}$ in the direction $v$. For $\kappa$ a small positive real number, to be chosen later on, we set
$$\Lambda_{z,v}=\{\lambda\in\cc,\ |\lambda|<3\kappa\tau(z,v,|\rho(z)|)\text{ and } z+\lambda v\in X\}.$$
The points $z+\lambda v$, $\lambda\in\Lambda_{z,v},$ are the points of $X$ which belong to $\Delta_{z,v}\left(3\kappa \tau(z,v,|\rho(z)|)\right)$, thus they all belong to $D$ provided $\kappa<\frac13$.\\
For $\lambda\in \Lambda_{z,v}$ let us define $g_{z,v}[\lambda]=g(z+\lambda v)$
and if $g_{z,v}[\lambda_1,\ldots,\lambda_k]$ is defined, let us set for $\lambda_1,\ldots, \lambda_k,\lambda_{k+1}$ belonging to $\Lambda_{z,v}$ and pairwise distinct
$$g_{z,v}[\lambda_1,\ldots,\lambda_{k+1}]=\frac{g_{z,v}[\lambda_1,\ldots, \lambda_k]-g_{z,v}[\lambda_2,\ldots,\lambda_{k+1}]}{\lambda_1-\lambda_{k+1}}.$$
Now consider the quantity
$$c_\infty(g)=\sup|g_{z,v}[\lambda_1,\ldots,\lambda_k]|  \tau(z,v,|\rho(z)|)^{k-1}$$
where the supremum is taken over all $z\in D,$ all $v\in\cc^n$ with $|v|=1$ and all $\lambda_1,\ldots,\lambda_k\in\Lambda_{z,v}$ pairwise distinct.
In Section \ref{section5}, we will prove that the finiteness of $c_\infty(g)$ implies the existence of a smooth extension $\tilde g$ which satisfies the hypothesis of Theorem \ref{th0}. We will then obtain the following theorem
\begin{theorem}\mlabel{th3}
In $\cc^2$, any function $g$ holomorphic on $X\cap D$  such that $c_\infty(g)$ is finite admits a holomorphic extension $G$ which belongs to $BMO(D)$ such that $\|G\|_{BMO(D)}$ is bounded up to a multiplicative uniform constant by $c_\infty(g)$.
\end{theorem}
Conversely, if we know that $g$ admits a bounded holomorphic extension $G$ on $D$ and if $\lambda_1,\lambda_2$ belong to $\Lambda_{z,v}$, Montel in \cite{Mon} proves that there exist a point $a$ in the unit disc of $\cc$ and $\mu$ in the segment $[\lambda_1,\lambda_2]$ such that 
$\frac{g_{z,v}(\lambda_1)-g_{z,v}(\lambda_2)}{\lambda_1-\lambda_2}$ can be written as $a \diffp{G}{v}(z+\mu v)$. But since $G$ is bounded, its derivative, and therefore
the divided difference $\frac{g_{z,v}(\lambda_1)-g_{z,v}(\lambda_2)}{\lambda_1-\lambda_2}$ as well,
are bounded by $\|G\|_{L^\infty(D)}$ times the inverse of the distance from $z+\mu v$ to the boundary of $D$ in the direction $v$, and this quantity is comparable to $\tau(z,v,|\rho (z)|)$. We will show in Section \ref{section5} that this necessary condition holds in fact in $\cc^n$, $n\geq 2$, and for more than two points $\lambda_1$ and $\lambda_2$, and so we will prove the following theorem 

\begin{theorem}\mlabel{th1}
In $\cc^n$, $n\geq 2$, if a function $g$ holomorphic on $X\cap D$ admits an  extension $G$ which is bounded and holomorphic on $D$ then $c_\infty(g)$ is finite.
\end{theorem}
In Section \ref{section5}, we will also study the case of $L^q$ extensions and, still using divided differences, we will give in $\cc^n$, $n\geq 2$, a necessary condition for a function $g$ holomorphic on $X\cap D$ to admit a holomorphic extension to $D$ which belong to $L^q(D)$. Then we will also prove that this condition is sufficient when $n=2$ (see Theorem \ref{th4}, \ref{th5} and \ref{th6} for precise statements). We will also see in Section \ref{section5}, Theorem \ref{th8} and \ref{th9}, that all these results can be generalized in a natural way to weak holomorphic functions in the sense of Remmert.

It should be noticed that a condition using divided differences was already used in \cite{Duquenoy} but that only varieties with singularities of order 2 were considered there. Here we have no restriction on the order of the singularities, and our condition uses all the divided differences of degree at most the orders of the singularities.

In Section \ref{section6}, we illustrate  these conditions by examples. Among other things,  when $D$ is the ball of center $(1,0)$ and radius $1$  and  $X=\{(z=(z_1,z_2)\in\cc^2, z_1^q=z_2^2\}$, with $q$ a positive odd integer, we will prove that any $g$ holomorphic and bounded on $X\cap D$ has a $L^2$-holomorphic extension on $D$ if and only if $q=1$ or $q=3$.
\par\bigskip
The article is organized as follows. In Section \ref{section2} we fix our notations and recall some results concerning the Berndtsson-Andersson kernel. In Section \ref{section3} we construct the new residue current adapted to our extension problem, and we prove Theorem \ref{th0} in Section \ref{secIII}. In Section \ref{section5} we prove Theorem \ref{th3} and \ref{th1} and we treat the case of $L^q$ holomorphic extensions.  We give examples of applications of our results in Section \ref{section6}.

\section{Notations and tools} \mlabel{section2}
As usually, when $BMO$ questions or estimates of integral kernels arise in this context, the Koranyi balls or McNeal polydiscs, their generalization for convex domains of finite type, naturally  appear (see \cite{AB,AC, DFF} for example). This will be of course the case in this article, but here  (and it seems to be the first time this happens) the Koranyi balls will appear directly in the construction of the residue current, and so in the construction of a good extension. These balls enable us to establish a connection between the geometric properties of the boundary of the domain and the geometric properties of the variety  (see Section \ref{section3}). The second classical tool we use is the Berndtsson-Andersson reproducing kernel which we also recall in this section.

\subsection{Notations}
Let us first fix our notation and adopt the following convention. We will often have estimates up to multiplicative constants. For readability convenience we introduce the following notation: We write $A\leqs B$ if there exists some constant $c>0$ such that $A\leq cB$. Each time we will mention from which parameters $c$ depends. We will write $A\eqs B$ if $A\leqs B$ and $B\leqs A$ both holds.

We write $X$ as $X=\{z,\ f(z)=0\}$ where $f$ is a holomorphic function defined in a neighborhood of $\overline D$. Without restriction we assumed that $f$ is minimal (see \cite{Chabat}, Theorem 3, paragraph 50).

\subsection{Koranyi balls}
We call the coordinates
system centered at $\zeta$ of basis $\eta_{\zeta}, v_{\zeta}$ the \ko coordinates
system at $\zeta$. We denote by $(z_1^*,z_2^*)$ the coordinates of a
point $z$ in the \ko coordinates system centered at $\zeta$. The \ko ball centered
in $\zeta$ of radius $r$ is the set $\kb r {\zeta}:=\{\zeta+\lambda\eta_{\zeta}+\mu
v_{\zeta},\ |\lambda|<r,\ |\mu|<r^{\frac12}\}$. These balls have the following properties~:
\begin{proposition}\mlabel{propII.0.1} There exists a neighborhood $\cal U$ of $bD$ and positive real numbers $\kappa$ and $c_1$ such that 
\begin{enumerate}[(i)]
 \item  for all $\zeta\in {\cal U}\cap D$, ${\cal P}_{4\kappa|\rho(\zeta)|}(\zeta)$ is included in $D$.
 \item  for all $\varepsilon>0$, all $\zeta,z\in {\cal U}$, $\p\varepsilon\zeta\cap\p\varepsilon z\neq
\emptyset$ implies $\p\varepsilon z\subset \p{c_1\varepsilon}\zeta$.
\item  for all $\varepsilon>0$ sufficiently small, all
$z\in {\cal U}$, all $\zeta\in \p{\varepsilon}z$ we have $|\rho(z)-\rho(\zeta )|\leq
c_1 \varepsilon$.
\item For all $\varepsilon >0$, unit vector $v\in\cc^n$, all $z\in{\cal U}$ and all $\zeta\in{\cal P}_\varepsilon (z)$, $\tau(z,v,\varepsilon )\eqs\tau(\zeta ,v,\varepsilon )$ uniformly with respect to $\varepsilon ,$ $z$ and $\zeta $.
\end{enumerate}
\end{proposition}
For $\cal U$ given by Proposition \ref{propII.0.1} and $z$ and $\zeta$ belonging to $\cal U$, we set $\delta(z,\zeta)=\inf\{\varepsilon>0, \zeta\in \p\varepsilon z\}$. Proposition \ref{propII.0.1} implies  that $\delta$ is a pseudo-distance in the
following sense: 
\begin{proposition}\label{propII.0.2}
For $\cal U$ and $c_1$ given by Proposition \ref{propII.0.1} and for all $z,\ \zeta$ and $\xi$ belonging to $\cal U$  we have
$$\frac1{c_1}\delta(\zeta,z)\leq \delta(z,\zeta)\leq c_1 \delta(\zeta,z)$$
 and
$$\delta(z,\zeta)\leq c_1(\delta(z,\xi)+\delta(\xi,\zeta))$$
\end{proposition}

\subsection{Berndtsson-Andersson reproducing kernel}
\mlabel{secII.0}
We now recall the definition of the Berndtsson-Andersson kernel of $D$ when $D$ is a  strictly convex domain of $\cc^2$.
We set $h_i(\zeta,z)=-\diffp{\rho}{\zeta_i}(\zeta)$, $h=\sum_{i=1,2} h_id\zeta_i$ and $\tilde{h}=\frac{1}{\rho}
h$. For a $(1,0)$-form $\beta (\zeta,z)=\sum_{i=1,2}\beta_i
d\zeta_i$ we set $\langle \beta (\zeta,z),\zeta-z\rangle = \sum_{i=1,2}
\beta_i(\zeta,z)(\zeta_i-z_i)$.
Then we define the Berndtsson-Andersson reproducing kernel by setting for
an arbitrary positive integer $N$, $n=1,2$ and all $\zeta,z\in D$~:
$$P^{N,n}(\zeta,z)=C_{N,n} \left(\frac{1}{1+\langle
\tilde{h}(\zeta,z),\zeta-z\rangle }\right)^{N+n}\left(\overline \partial
\tilde{h}\right)^n,$$
where $C_{N,n}\in \cc$ is a constant. We also set $P^{N,n}(\zeta,z)=0$ for all $z\in D$ and all $\zeta\notin D$. Then the following theorem holds (see
\cite{BA}):
\begin{theorem}
 For all $g\in \oo(D)\cap C^\infty(\overline D)$ we have
$$g(z)=\int_D g(\zeta)P^{N,2}(\zeta,z).$$
\end{theorem}
In the estimations of this kernel, we will need to write $h$ in the \ko coordinates at some point $\zeta_0$ belonging to $D$. We set  for $i=1,2$  $h_i^*=-\diffp{\rho}{\zeta^*_i}(\zeta)$. Then $h$ is equal to $\sum_{i=1,2}h_i^* d\zeta^*_i$ and satisfies the following proposition.
\begin{proposition}\mlabel{estiBA}
There exists a neighborhood $\cal U$ of $bD$ such that for all $\zeta\in D\cap {\cal U}$, all $\varepsilon>0$ sufficiently small
and all $z\in \p\varepsilon\zeta$ we have
\begin{enumerate}[(i)]
 \item $|\rho(\zeta)+\langle h(\zeta,z),\zeta-z\rangle|\geqs
\varepsilon+|\rho(\zeta)|+|\rho(z)|$,
 \item $|h^*_1(\zeta,z)|\leqs1$,
\item $|h^*_2(\zeta,z)|\leqs \varepsilon^{\frac{1}{2}}$,
\end{enumerate}
and there exists $c>0$ not depending from $\zeta$ nor from $\varepsilon$ such that
for all $z\in \p\varepsilon\zeta\setminus c\p\varepsilon\zeta$ we have
$$|\langle h(\zeta,z),\zeta-z\rangle|\geqs \varepsilon+|\rho(z)|+|\rho(\zeta)|,$$
uniformly \wrt $\zeta,z$ and $\varepsilon$.
\end{proposition}

\section{Construction of the extension operator}\mlabel{section3}
The holomorphic extension provided by  Theorem \ref{th0} will be given by a linear integral operator. Its definition is based upon the construction of Mazzilli
in \cite{Maz1} which uses Berndtsson-Andersson's reproducing kernel and  a current $T$ such that $fT=1$. The current $T$ relies on a family of currents $T_{\cal V}$, where ${\cal V}$ is an open subset of $D$, such that $fT_{\cal V}=1$.
Then using a locally finite covering $\left({\cal V}_j\right)_{j\in\nn}$ of $D$
and a partition of unity $\left(\chi_j\right)_{j\in\nn}$ associated with this
covering, Mazzilli glues together all the currents $T_{{\cal V}_j}$ and gets a current
$T=\sum_{j\in\nn} \chi_jT_{{\cal V}_j}$ such that $fT=1$. In \cite{Maz1}, the
only assumption on the covering $\left({\cal V}_j\right)_j$ is to be locally
finite.
\par\smallskip
In  order to get very fine estimates of the operator, instead of an ordinary locally finite covering, we will use a
covering of $D$ by Koranyi balls
$\left( {\cal P}_{\kappa |\rho(z_j)|}(z_j)\right)_{j\in\nn}$
which will be more suited to the geometry of $bD$ (see subsection \ref{maxcover}).
\par\smallskip
In \cite{Maz1}, the local current $T_{\cal V}$ is constructed using the Weierstrass polynomial $P_f$ of $f$ in the open set ${\cal V}$. This means that every
roots of $P_f$, or equivalently every sheets of $X$ intersecting  ${\cal V}$, are used. We will
modify the construction of $T_{\cal V}$  in order to use only the sheets of $X$ which
are meaningful for our purpose. In order to be able to choose the good sheets of $X$, we construct in  subsection \ref{secII.2} for $z_0$ near $bD$ a parametrization of $X$ in the \ko ball ${\cal P}_{\kappa|\rho(z_0)|}(z_0)$.
\par\smallskip
At last, we will have all the tools to define in subsection \ref{secII.3} the current $T$ such that $fT=1$ and the extension operator.

\subsection{Koranyi covering}\mlabel{maxcover}
In this subsection, for $\varepsilon_0>0$, we cover $D\setminus D_{-\varepsilon_0}$ with a family of \ko balls $\left({\cal P}_{\kappa|\rho(z_j)|}(z_j)\right)_{j\in\nn}$ where $\kappa$ is a positive small real number. This construction uses classical ideas of the theory of homogeneous spaces and is analogous to the construction of the covering of
\cite{BCD}.\\
Let $\varepsilon_0$, $\kappa$ and $c$ be positive real numbers sufficiently small. We construct a sequence of point of $D\setminus D_{\varepsilon _0}$ as follows.\\
Let $k$ be a non negative integer and choose $z_1^{(k)}$ in $bD_{-(1-c\kappa )^k\varepsilon _0}$ arbitrarily.\\
When $z_1^{(k)},\ldots, z_j^{(k)}$ are chosen, they are two possibilities. Either for all $z\in bD_{-(1-c\kappa )^k\varepsilon _0}$ there exists $i\leq j$ such that $\delta(z,z_i^{(k)})<c\kappa (1-c\kappa )^k\varepsilon _0$ and the process ends here or there exists $z\in bD_{-(1-c\kappa )^k\varepsilon _0}$ such that  for all $i\leq j$ we have $\delta(z,z_i^{(k)})\geq c\kappa (1-c\kappa )^k\varepsilon _0$ and we chose $z^{(k)}_{j+1}$ among these points. Since $D_{-(1-c\kappa )^k\varepsilon _0}$ is bounded, this process stops at some rank $n_k$.\\
We thus have constructed a sequence $(z_j^{(k)})_{k\in\nn, j\in\{1,\ldots,n_k\}}$ such that
\begin{enumerate}[(i)]
 \item \label{seqi} For all $k\in\nn$, and all $j\in\{1,\ldots, n_k\}$, $z_j^{(k)}$ belongs to $bD_{-(1-c\kappa )^k\varepsilon _0}$.
 \item \label{seqii} For all $k\in\nn$, all $i,j\in \{1,\ldots, n_k\}$, $i\neq j$, we have $\delta(z_i^{(k)},z_j^{(k)})\geq c\kappa (1-c\kappa )^k\varepsilon _0$.
 \item \label{seqiii} For all $k\in\nn$, all $z\in bD_{-(1-c\kappa )^k\varepsilon _0}$, there exists $j\in\{1,\ldots, n_k\}$ such that $\delta(z,z_j^{(k)})<c\kappa (1-c\kappa )^k\varepsilon _0$.
\end{enumerate}
For such sequences, we prove the following proposition.
\begin{proposition}\mlabel{propmax}
For $\kappa >0$ and $c>0$ small enough, let $\left(z_j^{(k)}\right)_{k\in\nn,j\in\{1,\ldots, n_k\}}$ be a sequence which satisfies (\ref{seqi}), (\ref{seqii}) and (\ref{seqiii}). Then
\begin{enumerate}[(a)]
 \item \label{propmax1} $D\setminus D_{\varepsilon _0}$ is included in $ \cup_{k=0}^{+\infty} \cup_{j=1}^{n_k} {\cal P}_{\kappa |\rho (z_j^{(k)})|}\left(z_{j}^{(k)}\right)$,
 \item \label{propmax2} there exists $M\in\nn$ such that for $z\in D\setminus D_{-\varepsilon _0}$, ${\cal P}_{4\kappa |\rho (z)|}(z)$ intersect at most $M$ Koranyi balls ${\cal P}_{ 4\kappa |\rho (z_j^{(k)})|}\left(z_{j}^{(k)}\right)$.
\end{enumerate}
\end{proposition}
\pr 
We first prove that (\ref{propmax1}) holds. For $z\in D\setminus D_{\varepsilon _0}$ let $k\in\nn$ be such that 
$$(1-c\kappa )^{k+1}\varepsilon _0<|\rho (z)|<(1-c\kappa )^k\varepsilon _0$$
and let $\lambda\in\cc$ be such that $\zeta=z+\lambda \eta_z$ belong to $bD_{-(1-c\kappa)^k\varepsilon_0}$. 
On the one hand the assumption $(\ref{seqiii})$ implies that there exists $j\in\{1,\ldots, n_k\}$ such that $\delta\left(\zeta ,z_j^{(k)}\right)\leq c\kappa (1-c\kappa )^k\varepsilon _0$. On the other one hand we have $|\lambda|=\delta(z,\zeta)\leq C c\kappa(1-c\kappa)^k\varepsilon_0$ where $C$ does not depend from $z$ nor from  $\zeta $. These two inequalities yield
\begin{eqnarray*}
 \delta\left(z,z^{(k)}_j\right)&\leq& c_1(\delta(z,\zeta )+c_1\delta(\zeta ,z^{(k)}_j)\\
&\leq& \kappa cc_1(1-c\kappa )^k\varepsilon _0 (C\kappa +1)\\
&\leq& \kappa |\rho \left(z_j^{(k)}\right)|
\end{eqnarray*}
provided $c$ is small enough. Therefore $z$ belongs to ${\cal P}_{\kappa |\rho (z_j^{(k)})|}(z_j^{(k)})$ and (\ref{propmax1}) holds.\\
We now prove (\ref{propmax2}). Let $z$ be a point of $D\setminus D_{\varepsilon _0}$. For all $\zeta \in{\cal P}_{4\kappa |\rho(z)|}(z)$, if $\kappa $ is small enough, proposition \ref{propII.0.1} yields
$$\frac12 |\rho (z)|\leq |\rho (\zeta )|\leq 2|\rho (z)|.$$
The same inequalities hold for all $z^{(k)}_j$ and all $\zeta\in {\cal P}_{4\kappa |\rho (z_j^{(k)})|}(z_j^{(k)})$. 
Thus if ${\cal P}_{4\kappa |\rho (z_j^{(k)})|}(z_j^{(k)})\cap {\cal P}_{\kappa |\rho (z)|}(z)\neq \emptyset$ we have
$$\frac14|\rho (z)| \leq (1-c\kappa )^k\leq 4 |\rho (z)|.$$
Therefore $k$ can take at most $\frac{4\ln 2}{|\ln(1-c\kappa )|}$ values.\\
For such a $k$, we set $I_k=\left\{j\in\{1,\ldots,n_k\},\ {\cal P}_{4\kappa |\rho (z_j^{(k)})|}(z_j^{(k)})\cap {\cal P}_{4\kappa |\rho (z)|}(z)\neq \emptyset\right\}$. Assertion (\ref{propmax2}) will be proved provided we show that $\#I_k$, the cardinal of $I_k$, is bounded uniformly with respect to $k$ and $z$.\\
We denote by $\sigma$ the area measure on $bD_{-(1-c\kappa )^k\varepsilon _0}$. Since for all $i,j\in I_k$ distinct we have
$ \delta\left(z_{i}^{(k)},z_j^{(k)}\right)\geq c\kappa (1-c\kappa )^k\varepsilon _0$, provided $c$ is small enough,  we have
\begin{eqnarray*}
\lefteqn{\sigma \left(\cup_{j\in I_k}{\cal P}_{4\kappa \left|\rho \left(z_j^{(k)}\right)\right|}\left(z_j^{(k)}\right)\cap bD_{-(1-c\kappa )^k\varepsilon _0}\right)}\\
& \geq&\sigma \left(\cup_{j\in I_k}{\cal P}_{\frac c{c_1} \kappa(1-c\kappa )^k\varepsilon _0}\left(z_j^{(k)}\right)\cap bD_{-(1-c\kappa )^k\varepsilon _0}\right)\\
& \geq& \# I_k \left(\frac c{c_1} \kappa (1-c\kappa )^k\varepsilon _0\right)^n.
\end{eqnarray*}
Now we look for an upper bound of $\sigma \left(\cup_{j\in I_k}{\cal P}_{4\kappa |\rho (z_j^{(k)})|}(z_j^{(k)})\cap bD_{-(1-c\kappa )^k\varepsilon _0}\right)$.
We  fix $j_0\in I_k$. For all $j\in I_k,$ since ${\cal P}_{4\kappa |\rho (z_j^{(k)})|}(z_j^{(k)})\cap {\cal P}_{4\kappa |\rho (z)|}(z)\neq \emptyset$ and ${\cal P}_{4\kappa |\rho (z_{j_0}^{(k)})|}(z_{j_0}^{(k)})\cap {\cal P}_{4\kappa |\rho (z)|}(z)\neq \emptyset$, we have 
\begin{eqnarray*}
 \delta\left(z_{j_0}^{(k)},z_j^{(k)}\right)&\leqs & \delta\left(z_{j_0}^{(k)},z\right)
 +\delta\left(z,z_j^{(k)}\right)\\
&\leqs& 4\kappa  \left(\left|\rho \left(z_{j_0}^{(k)}\right)\right|+\left|\rho \left(z_{j}^{(k)}\right)\right|\right)\\
&\leqs & \kappa (1-c\kappa )^k\varepsilon _0
\end{eqnarray*}
uniformly with respect to $k$, $j$ and $j_0$. Thus there exists $K$ not depending from $z$, $j$, $j_0$ nor on $k$ such that 
${\cal P}_{4\kappa |\rho (z_j^{(k)})|}(z_j^{(k)})
\subset {\cal P}_{\kappa K |\rho (z_{j_0}^{(k)})|}(z_{j_0}^{(k)})$. Therefore 
\begin{eqnarray*}
\sigma \left(\cup_{j\in I_k}{\cal P}_{4\kappa |\rho (z_j^{(k)})|}(z_j^{(k)})\cap bD_{-(1-c\kappa )^k\varepsilon _0}\right) 
&\leq&
\sigma \left({\cal P}_{4K\kappa |\rho (z_{j_0}^{(k)})|}(z_{j_0}^{(k)})\cap bD_{-(1-c\kappa )^k\varepsilon _0}\right)\\
&\leqs& \left(K\kappa (1-c\kappa )\varepsilon _0\right)^n 
\end{eqnarray*}
which yields $\#I_k\leqs c^{-n}$.\qed\\
The covering property (\ref{propmax1}) allows us to settle the following definition 
\begin{definition}
 Let $\cal U$ be any subset of $\cc^n$. If the sequence $(z_j)_{j\in\nn}$ can renumbered such that (\ref{seqi}), (\ref{seqii}) are satisfied and such that (\ref{seqiii}) holds true for all $z\in {\cal U}\cap (D\setminus D_{-\varepsilon_0})$, the family $\left({\cal P}_{\kappa |\rho (z_j)|}(z_j)\right)_{j\in\nn}$ will be called a $\kappa$-covering of ${\cal U}\cap (D\setminus D_{-\varepsilon_0})$.
\end{definition}
\subsection{A family of parametrizations}\mlabel{secII.2}
In order to construct the current we need to define our extension operator, we will need some kind of parametrization for $X$ over ${\cal P}_{\kappa|\rho(z_0)|}(z_0)$ when $z_0$ is near the boundary of the domain and when ${\cal P}_{\kappa|\rho(z_0)|}(z_0)\cap X\neq \emptyset$. Moreover, we will need some uniform estimates for this parametrization. Of course if we are near a regular point of $X$, such parametrizations do exist but the situation is more delicate when we are near a singularity of $X$. Given a point $z_0$ near a singularity $\zeta_0$ of $X$ which belongs to $bD$, we denote by $(\zeta_{0,1}^*,\zeta^*_{0,2})$ the coordinates of $\zeta_0$ is the \ko coordinates at $z_0$. We denote by $\Delta$ the unit of $\cc$ and by $\Delta_z(r)$ the disc of $\cc$ centered at $z$ of radius $r$. Our goal in this subsection is  to prove the following propositions:
\begin{proposition}\label{new_prop}
There exists $\kappa >0$ sufficiently small and not depending on $z_0$ such that if $X\cap {\cal P}_{\kappa|\rho(z_0)|}(z_0)\neq\emptyset,$ then $|\zeta^*_{0,1}|\geq 2\kappa|\rho(z_0)|$.
\end{proposition}
\begin{proposition}\mlabel{propII.2.1}
 There exist $\kappa$ and $r$ positive real numbers sufficiently small, a positive integer $p_0$ and a neighborhood $\cal U$ of $\zeta_0$ such that for all $z_0\in {\cal U}$, if $|\zeta^*_{0,1}|\geq\kappa|\rho(z_0)|$ then there exist 
$\alpha_1^*,\ldots, \alpha^*_{p_0}$ holomorphic functions in $\Delta_0(2\kappa|\rho(z_0)|)$ which satisfy 
\begin{enumerate}[(i)]
 \item \mlabel{propII.2.1.2} $\alpha_j^*$ and
$\diffp{\alpha^*_j}{z^*_1}$ are bounded on $\Delta_0(2\kappa|\rho(z_0)|)$ uniformly \wrt $z_0.$
\item  \mlabel{propar3}if there exists $j$ and $z^*_1$ such that
$(z^*_1,\alpha_j^*(z^*_1))$ belong to ${\cal P}_{2\kappa|\rho(z_0)|}(z_0)$ then for all $\zeta_1^*\in \Delta_0(2\kappa|\rho(z_0)|)$ we have $|\alpha_j^*(\zeta_1^*)|\leq
\left(3\kappa|\rho(z_0)|\right)^{\frac12}.$
\item \mlabel{propII.2.1.4} There exists $u_0$ holomorphic in $\Delta_{z_0}(r)^2$ such that $|u_0|\eqs 1$ uniformly with respect to $z_0$ and $f(\zeta)=u_0(\zeta) \prod_{i=1}^{p_0}(\zeta^*_2-\alpha_{i}^*(\zeta^*_1))$ for all $\zeta\in{\cal P}_{2\kappa|\rho(z_0)|}(z_0)$.
\end{enumerate}
\end{proposition}
The proofs of this proposition will relies on the following two lemmas. 
\begin{lemma}\mlabel{propweierstrass}
 Let $(A,d)$ be a metric space, $\alpha_0\in A$ and $(f_\alpha)_{\alpha\in A}$ a
family of holomorphic function on $\Delta^2$ such that
\begin{itemize}
 \item[-] $(f_\alpha)_{\alpha\in A}$ converges uniformly to $f_{\alpha_0}$ when
$\alpha$ tends to $\alpha_0$,
 \item[-] $f_{\alpha_0}(0,\cdot)\neq 0$ and $f_{\alpha_0}(0)=0$.
\end{itemize}
Then there exist positive real numbers $r_1,r_2,\eta>0$, a positive integer $p$ such that, for all $\alpha\in A$ with $d(\alpha,\alpha_0)<\eta$, there exist $p$ functions $a_1^{(\alpha)},\ldots,a_p^{(\alpha)}$ holomorphic on $\Delta_0(r_1)$ and a function $u_\alpha$ holomorphic in $\Delta_0(r_1)\times \Delta_0(r_2)$  which satisfy
\begin{enumerate}[(i)]
 \item $f_\alpha(z)=u_\alpha(z) \left(z_2^p+a_1^{(\alpha)}(z_1)
z_2^{p-1}+\ldots+ a_p^{(\alpha)}(z_1)\right)$,
 \item $|u_\alpha(z)|\eqs 1$ for all $z\in \Delta_0(r_1)\times \Delta_0(r_2)$
uniformly \wrt $z$ and $\alpha$.
\end{enumerate}
\end{lemma}
\pr 
We first want to apply Rouch\'e's theorem to
$f_{\alpha}(z_1,\cdot)-f_{\alpha_0}(0,\cdot)$, $z_1$ fixed in $\Delta_0(r_1)$ where
$r_1>0$ is to be chosen in a moment.\\
Since $f_{\alpha_0}(0,\cdot)$ is not identically zero, there exists $r_2>0$ such that $f_{\alpha_0}(0,z_2)\neq 0$ for
all $z_2\in\Delta_0(r_2)\setminus\{0\}$. We denote by $a$ the positive real number $a=\inf_{|z_2|=r_2}|f_{\alpha_0}(0,z_2)|$ and by $p$ the order the root $0$ of $f_{\alpha_0}(0,\cdot)$.\\
Since $(f_\alpha )_\alpha $ converges uniformly to $f_{\alpha _0}$ on $\Delta_0(1)$, there exists $\eta>0$ such that for all $\alpha\in A$,
$d(\alpha_0,\alpha)<\eta$, all $z\in\Delta_0(1)^2$ the following inequality holds:
$\sup_{z\in \Delta_0(1)^2}|f_\alpha(z)-f_{\alpha_0}(z)|<\frac a4$.\\
By Cauchy's inequalities, there exists $r_1>0$ such that  for all $z\in
\Delta_0(r_1)\times \Delta_0(r_2)$ we have
$|f_{\alpha_0}(z_1,z_2)-f_{\alpha_0}(0,z_2)|<\frac a4$.\\
Thus $|f_\alpha(z_1,z_2)-f_{\alpha_0}(0,z_2)|\leq |f_{\alpha_0}(0,z_2)|$ and by
Rouch\'e's theorem, $f_\alpha(z_1,\cdot)$ has exactly $p$ zeros in $\Delta_0(r_2)$
for all $z_1$ fixed in $\Delta_0(r_1)$. Therefore by the Weierstrass preparation theorem there exist $p$ functions 
$a_1^{(\alpha)}, \ldots, a_p^{(\alpha)}$ holomorphic on $\Delta_0(r_1)$ and a function $u_\alpha$ holomorphic on $\Delta_0(r_1)\times \Delta_0(r_2)$ zero free such that
$$f_\alpha(z)=u_\alpha(z) \left(z_2^p+a^{(\alpha)}_1(z_1)
z_2^{p-1}+\ldots+a^{(\alpha)}_p(z_1)\right).$$
We set $P_\alpha(z_1,z_2)=z_2^p+a^{(\alpha)}_1(z_1) z_2^{p-1}+\ldots+a^{(\alpha)}_p(z_1)$.
To end the proof of the lemma we have to prove that $1\leqs| u_\alpha|\leqs
1$. We prove the lower uniform boundedness.\\
For all $z_1\in \Delta_0(r_1)$, $\frac1{u_\alpha(z_1,\cdot)}$ is holomorphic and
$$\frac{1}{|u_\alpha(z_1,z_2)|}\leq
\max_{|\zeta_2|=r_2}\left|\frac{P_\alpha(z_1,\zeta_2)}{f_\alpha(z_1,\zeta_2)}
\right|.$$
On the one hand, for all $\alpha\in A$ such that $d(\alpha,\alpha_0)<\eta$, all
$(z_1,z_2)\in \Delta_0(r_1)\times b\Delta_0(r_2)$ we have
\begin{eqnarray*}
|f_\alpha(z)|&\geq&
|f_{\alpha_0}(0,z_2)|-|f_{\alpha_0}(z)-f_{\alpha_0}(0,z_2)|-|f_\alpha(z)-f_{
\alpha_0}(z)|\\
&\geq& a-\frac a4-\frac a4=\frac a2.
\end{eqnarray*}
On the other one hand, since $(f_\alpha)_{\alpha\in A}$ converges uniformly to
$f_{\alpha_0}$ when $\alpha$ tends to $\alpha_0$ and since $f_\alpha(z)$ is
uniformly bounded away from $0$ for $(z_1,z_2)\in\Delta_0(r_1)\times b\Delta_0(r_2)$,
$(a_j^{(\alpha)})_{\alpha\in A}$ converge uniformly to $a_j^{(\alpha_0)}$ for
all $j$ when $\alpha$ tends to $\alpha_0$. This implies that $(P_\alpha)_{\alpha\in A}$
converges uniformly to $P_{\alpha_0}$ and therefore $\sup_{\Delta_0(r_1)\times
\Delta_0(r_2)} |P_\alpha|$ is uniformly bounded for $\alpha$ near $\alpha_0$.\\
This yields $|u_\alpha(z)|\geqs 1$ uniformly \wrt $z\in\Delta_0(r_1)\times \Delta_0(r_2)$ and
$\alpha\in A$ such that $d(\alpha,\alpha_0)<\eta$. The upper boundedness can be
proved in the same way.\qed
\begin{lemma}\mlabel{lemII.2.2}
Let $\zeta_0\in bD$ be a singularity of $X$, let $z_0\in D$ be a point near enough $\zeta_0$. 
There exist $r>0$ not depending from
$z_0$ and a parametric  representation of $X$ in the \ko coordinates system
centered at $z_0$ of the form $({t^*}^p+\zeta^*_{0,1}, \phi(t^*)+\zeta^*_{0,2})$,
 such that $|\phi^*(t^*)|\leqs \left|t^*\right|^p$, $t^*\in \Delta_0(r)$, uniformly
with respect to $z_0$. 
\end{lemma}
\pr Without restriction we assume that $\zeta _0$ is the origin of $\cc^2$. Maybe after a unitary linear change of coordinates, there exists $r_0>0$, $p,q\in\nn$, $q>p>1$, and $u$ holomorphic and bounded on $\Delta_0(r_0)$, $u(0)\neq 0$ such
that $\phi:t\mapsto (t^p,t^qu(t))$ is a parametric representation of $X$ over
$\Delta_0(r_0)$.\\
We consider $z_0$ such that $|\zeta_0-z_0|<r_0$ and we denote by $(\alpha,\beta)$ the coordinates of $\eta_{z_0}$ and by
$(-\overline{\beta},\overline{\alpha})$ the coordinates of $v_{z_0}$.
In the \ko coordinates centered at $z_0$, $X$ is parametrized by $t\mapsto
(\overline{\alpha}t^p+\overline{\beta}t^qu(t)+\zeta^*_{0,1}, -\beta t^p+\alpha
t^qu(t)+\zeta^*_{0,2})$.\\
Let $(\alpha_0,\beta_0)$ denotes the coordinates of $\eta_{\zeta_0}$. The transversality
hypothesis implies that $\alpha_0\neq 0$ so there exists $r_1>0$ and a $p$-th determination of the root $\phi_1$ in $\Delta_{\overline{\alpha_0}}(r_1).$
If $r_0>0$ is sufficiently small, ${\alpha}$ belongs to
$\Delta_{\alpha_0}(r_1)$ and 
$\overline{\alpha}t^p+\overline{\beta}t^q u(t)=(\phi_1(\overline{\alpha})t)^p
\left(1+\frac{\overline{\beta}}{\overline{\alpha}}t^{q-p}u(t)\right).$\\
Since $q>p$, there exists $r_2\in]0,r_1[$ such that for all $t\in \Delta_0(r_2),$ all
$\beta\in \Delta_{\beta_0}(r_2)$ and all $\alpha\in \Delta_{\alpha_0}(r_2)$, we have
$\left|1+\frac{\overline{\beta}}{\overline{\alpha}}t^{q-p}u(t)  \right|\geq
\frac12$ and so there exists $\phi_2$ holomorphic for $t\in \Delta_0(r_2)$,
$C^\infty$-smooth for $\alpha\in \Delta_{\alpha_0}(r_2)$ and $\beta\in\Delta_{\beta_0}(r_2)$
such that
$\phi_2(t,\alpha,\beta)^p=1+\frac{\overline{\beta}}{\overline{\alpha}}t^{q-p}
u(t)$.\\
We apply the implicit functions theorem to $\Psi:
(t,t^*,\alpha,\beta)\mapsto t^*-\phi_1(\overline{\alpha})\phi_2(t,\alpha,\beta)
t$. Since $\Psi(0,0,\alpha_0,\beta_0)=0$ and
$\diffp{\Psi}{t}(0,0,\alpha_0,\beta_0)\neq 0$, there exist $r>0$ and
${\tilde \psi}:\Delta_0(r)\times \Delta_{\alpha_0}(r)\times \Delta_{\beta_0}(r)\to V(0)$, $V(0)$
neighborhood of $0\in \cc$ such that $\tilde \psi$ is holomorphic in $t$, and
$C^\infty$-smooth in $\alpha$ and $\beta$ such that ${t^*}^p=\overline{\alpha}
t^p+\overline{\beta} t^q u(t)$ if and only if $t=\tilde \psi(t^*,\alpha,\beta)$.\\
We now end the proof of the lemma by setting
$$\phi^*(t^*)=-\beta\tilde \psi(t^*,\alpha,\beta)^p+\alpha\tilde\psi(t^*,\alpha,
\beta)^qu\left(\tilde\psi(t^*,\alpha,\beta)\right).$$\qed\\[10pt]
{\it Proof of proposition \ref{new_prop}:}
We first choose $\kappa >0$ such that $2\kappa|\rho(z_0)|\leq r$, $r$ given by lemma
\ref{lemII.2.2} and we write $\zeta\in X\cap \p{\kappa|\rho(z_0)|}{z_0}$ as
$\zeta=\left({t^*}^{p_0}+\zeta_{0,1}^*,\phi^*(t^*)+\zeta^*_{0,2}\right)$ for some $t^*$ belonging to $\Delta_0(r)$. Now, if we assume
that 
$\left|\zeta^*_{0,1}  \right|< 2\kappa|\rho(z_0)|$ we get $|\zeta_1^*-\zeta^*_{0,1}|\leq
3\kappa|\rho(z_0)|$ and therefore $|t^*|\leq (3\kappa|\rho(z_0)|)^{\frac1{p_0}}$. This yields  
\begin{eqnarray*}
|\zeta^*_{0,2}|&\leq& |\zeta^*_{0,2}-\zeta^*_2|+|\zeta_2^*|\\
&\leq& |\phi^*(t^*)|+|\zeta_2^*|\\
&\leqs& \kappa|\rho(z_0)|+(\kappa|\rho(z_0)|)^{\frac12}\\
&\leqs& (\kappa|\rho(z_0)|)^{\frac12}
\end{eqnarray*}
uniformly with respect to $z_0$. Thus there exists $K>0$ not depending from $z_0$
nor from $\kappa $ such that $\zeta_0$ belongs to $\p{\kappa K|\rho(z_0)|}{z_0}$. Moreover, if $\kappa $ is
chosen sufficiently small, for all $\xi\in \p{\kappa K|\rho(z_0)|}{z_0}$ Proposition \ref{propII.0.1} gives 
$|\rho(\xi)|\geq \frac12 |\rho(z_0)|$. This gives a contradiction
because $|\rho(\zeta_0)|=0<|\rho(z_0)|$ whereas $\zeta_0$ belongs to $\p{\kappa K|\rho(z_0)|}{z_0}$. Therefore
we can choose $\kappa >0$ not depending from $z_0$ such that $\left|\zeta^*_{0,1} 
\right|\geq 2\kappa|\rho(z_0)|$.\qed\\[10pt]
{\it Proof of proposition \ref{propII.2.1}:} Let $p_0$ be the multiplicity of the singularity $\zeta_0$ of $X$ and 
let $\psi$
be a $p_0$-th determination of the root holomorphic in $\Delta_{\zeta^*_{0,1}}(2\kappa|\rho(z_0)|)$. We set
$\alpha_j^*(z_1^*)=\phi^*\left(\psi(z^*_1-\zeta^*_{0,1})
e^{\frac{2i\pi}{p_0}j}\right)+\zeta^*_{0,2}$, $j=1,\ldots, p_0$. For all $j$,
$\alpha_j^*$ is holomorphic on $\Delta_0(2\kappa|\rho(z_0)|)$ and  is
uniformly bounded on $\Delta_0(2\kappa|\rho(z_0)|)$. We have
$$\diffp{\alpha_j^*}{z_1^*}(z_1^*)=\psi'(z_1^*-\zeta^*_{0,1})
\diffp{\phi^*}{t^*}\left(\psi(z^*_1-\zeta^*_{0,1})e^{\frac{2i\pi}{p_0}j}\right)e^{
\frac{2i\pi}{p_0}j}.$$
Since $|\phi^*(t^*)|\leqs |t^*|^p$ this yields
$\left|\diffp{\alpha_j^*}{z_1^*}(z_1^*)  \right|\leqs 1$ which proves (\ref{propII.2.1.2}).
\par\medskip
We now prove that (\ref{propar3}) holds. We denote by $K$ a uniform bound of the
derivative of $\alpha_j^*$. If $z^*_1\in \Delta_0(2\kappa|\rho(z_0)|)$ is such that $|\alpha_j^*(z_1^*)|\leq
\left(2\kappa|\rho(z_0)|\right)^{\frac{1}{2}}$, we have for all $\zeta_1^*\in \Delta(2\kappa |\rho (z_0)|)$:
\begin{eqnarray*}
|\alpha_j^*(\zeta_1^*)|&\leq&
|\alpha_j^*(z_1^*)|+\left|\alpha_j^*(z_1^*)-\alpha_j^*(\zeta_1^*)\right|\\
&\leq&(2\kappa|\rho(z_0)|)^{\frac12}+K|\zeta_1^*-z_1^*|\\
&\leq& (2\kappa|\rho(z_0)|)^{\frac12}+4K\kappa|\rho(z_0)|.
\end{eqnarray*}
Therefore choosing again $\kappa $ small enough, uniformly with respect to $z_0$, we get 
$|\alpha_j^*(\zeta_1^*)|\leq \left(3\kappa|\rho(z_0)|\right)^{\frac{1}{2}}$.
\par\medskip
Only (\ref{propII.2.1.4}) is left to be shown. 
For $z$ near $\zeta_0$ we set $f_z(\lambda,\mu)=f(\zeta_0+\lambda \eta_{z}+\mu v_{z})$ and we apply
Lemma \ref{propweierstrass} to the family $(f_z)_z$ which gives
$u_0$ and $P_0$ such that $f_{z_0}=u_0P_0$ where $|u_0|\eqs 1$ uniformly with respect to $z_0$ and where $P_0(\lambda\eta_{z_0}+\mu v_{z_0})$ is a polynomial of the variable $\mu$ with coefficients holomorphic \wrt $\lambda$. We have $f_{z_0} (z_0-\zeta_0+\zeta_1^*\eta_{z_0}+\alpha_i^*(\zeta_1^*)v_{z_0})=0$ for all $i$ so
for all $\zeta$ such that $|\zeta^*_1|< 2\kappa|\rho(z_0)|$
$$P_0(\zeta_1^*-\zeta_{0,1}^*,\zeta_2^*-\zeta_{0,2}^*)=\prod_{i=1}^{p_0} (\zeta^*_2-\alpha_i^*(\zeta_1^*)).$$
\qed

\subsection{Definition of the operator}\mlabel{secII.3}
We now come to the definition of the current $T$ such that $fT=1$ and of the extension operator. Our construction is a refinement
of \cite{Maz1}. We choose a positive real number $\kappa$ so that  Propositions \ref{propmax} and \ref{propII.2.1}
 hold true for such a $\kappa $ and such that Proposition \ref{propII.0.1} implies that $2\rho(z)\leq \rho(\zeta)\leq\frac12\rho(z)$ for all $z\in D$ near $bD$.\\
For $\varepsilon _0>0$ and $z_0\in \overline{D_{-\varepsilon_0}}$, that is when $z_0$ is far from the boundary, we do not modify the construction except that we require that ${\cal U}_0$ is included in $D_{-\frac{\varepsilon_0}2}$. We get a covering ${\cal U}_{-m},\ldots, {\cal U}_{-1}$ of $\overline{D_{-\varepsilon_0}}$ and the corresponding currents $T_{-m},\ldots, T_{-1}$ such that $fT_j=1$ on ${\cal U}_j$ for all $j=-m,\ldots, -1$.\\
Near the boundary, we have to be more precise and we use a $\kappa $-covering $\left({\cal P}_{\kappa|\rho(z_j)|}(z_j)\right)_{j\in\nn}$ of $D\cap D_{-\varepsilon _0}$  constructed in Section \ref{maxcover}.
In the \ko coordinates centered at $z_j$, the fiber of $X$ above $(z^*_1,0)\in
\p{\kappa|\rho(z_j)|}{z_j}$ is given by $\{(z_1^*,\alpha^*_i(z_1^*)),\ i=1,\ldots, p_j\}$
where $p_j$ and $\alpha^*_1,\ldots, \alpha_{p_j}^*$ are given by Proposition
\ref{propII.2.1}. In \cite{Maz1}, Mazzilli actually considered the Weierstrass
polynomial in a neighborhood of $z_j$ but this neighborhood may be smaller than
$\p{\kappa|\rho(z_j)|}{z_j}$ or the Weierstrass polynomial may include all the
$\alpha_i^*$. However, in order to make a good link between the geometry of the boundary of $D$ and $X$, we need to have a polynomial in all $\pk{j}$ and we have to take into account only the sheets of $X$ which intersect $\pk{j}$ or equivalently the $\alpha_i^*$ such that for some $z_1^*\in \Delta_0(\kappa|\rho (z_j)|)$, the point  $z_j+z^*_1\eta_{z_j+}\alpha^*_i(z^*_1)v_{z_j}$ belongs to $\pk j$.
So we put $I_j\hskip -1.5pt =\hskip -1.5pt\left\{i, \exists z_1^*\in \Delta_0(\kappa
|\rho(z_j)|) \text{ such that } |\alpha_i^*(z^*_1)|\leq
(2\kappa|\rho(z_j)|)^{\frac12}\right\}$, $q_j=\#I_j$, the cardinal of $I_j$, and for any $C^\infty$-smooth
$(2,2)$-form $\phi$ compactly supported in $\pk j$ we set
$$\tilde{T}_j[\phi]=\int_{{\cal P}_{\kappa|\rho(z_j)|}(z_j)} \frac{\prod_{i\in
I_j}\overline{\zeta_2^*-\alpha_i^*(\zeta_1^*)}}{f(\zeta)}
\diffp{^{q_j}\phi}{\overline{\zeta^*_2}^{q_j}}(\zeta).$$
As in \cite{Maz1}, integrating by parts $q_j$-times gives $f\tilde T_j=c_j$ where $|c_j|=q_j!$.

Now let $\left(\chi_j\right)_{j\geq -m}$ be a partition of unity subordinated to
the covering ${\cal U}_{-m},\ldots, {\cal U}_{-1}$, $\left(\pk j\right)_{j\in\nn}$  of $D$. We assume that $\chi_j$ has been
chosen so that
$\left|\diffp{^{\alpha+\overline\alpha+\beta+\overline\beta}\chi_j}{{\zeta^*_1}
^\alpha \partial \overline{\zeta^*_1}^{\overline\alpha}
\partial{\zeta^*_2}^\beta
\partial\overline{\zeta^*_2}^{\overline\beta}}(\zeta)\right|\leqs
\frac{1}{|\rho(z_j)|^{\alpha+\overline\alpha+\frac{\beta+\overline\beta}{2}}}$ for
all $j\in\nn$, $\zeta\in\pk j$, $\alpha,\beta, \overline\alpha,\overline\beta\in\nn,$ uniformly with respect to $z_j$ and $\zeta$.
We set as in \cite{Maz1}: $T_j=\frac1{c_j} \tilde T_j$ for $j\in\nn$ and $T=\sum_{j=-m}^\infty \chi_j T_j$.

Therefore we have $fT=1$ on $D$. Moreover, since $T$ is supported in $\overline D$ which is compact, $T$ is of finite order
(see \cite{Sch}) and we can apply $T$ to smooth forms vanishing to a
sufficient order $l$ on $bD$. Therefore if the function $\tilde g$ is such that $|\rho|^N \tilde g$ belongs to $C^l(\overline D)$,  we can apply $T$ to $\tilde gP^{N,2}$. This gives us the integer $l$ of Theorem \ref{th0}.
\par\medskip
Let $b(\zeta,z)=\sum_{j=1,2}b_j(\zeta,z)d\zeta_j$ be the holomorphic $(1,0)$-form defined by $b_j(\zeta,z)=\int_0^1\diffp{f}{\zeta_j}(\zeta+t(z-\zeta))dt$ so that for all $z$ and $\zeta$ we have
$f(z)-f(\zeta)=\sum_{i=1,2}b_i(\zeta,z)(z_i-\zeta_j).$ 
Let $g$ be a holomorphic function admitting a smooth extension $\tilde g$ which satisfies the assumptions of Theorem \ref{th0}. Following the construction of \cite{Maz1}, we define the extension $E_N(g)$ of $g$ by setting 
$${E_N}[g](z)=C_1 \overline\partial T[\tilde gb(\cdot,z)\wedge P^{N,1}(\cdot,z)],\qquad\forall z\in D,$$
where $C_1$ is a suitable constant (see \cite{Maz1}). We have to check that $E_N(g)$ is indeed an extension of $g$.

We have the two following facts :\\
{\it Fact 1~:}\label{fact1} Mazzilli proved in \cite{Maz1} that if $\tilde g$ is holomorphic on $D$ and of class $C^l$ on $\overline D$ then 
${E_N}\tilde g=\tilde g$ on $X\cap D$.\\[2pt]
{\it Fact 2~:} \label{fact2}  We have $E_N\tilde g_1=E_N\tilde g_2$ when $\tilde g_1$ and $\tilde g_2$ are any smooth functions such that 
$\diffp{^{\alpha+\beta}\tilde g_1}{\overline{\zeta_1^*}^\alpha\partial \overline{\zeta_2^*}^\beta}=\diffp{^{\alpha+\beta}\tilde g_2}{\overline{\zeta_1^*}^\alpha\partial \overline{\zeta_2^*}^\beta}$ on $X\cap D$ for all integers $\alpha ,\beta $ with $\alpha +\beta \leq k$, where $k$ is the supremum of the orders of the singularities of $X$. Indeed, since $f$ is assumed to be minimal, using Theorem I, paragraph 11.2 and the theorem of paragraph 14.2 of \cite{Tsi}, for any function $\tilde g$ we can write $E_N \tilde g$ as a sum of integrals over $X\cap D$ where only the derivatives $\diffp{^{\alpha+\beta}\tilde gP^{N,1}}{\overline{\zeta_1^*}^\alpha\partial \overline{\zeta_2^*}^\beta}$ with $\alpha+\beta\leq k$. Applying this formula to $\tilde g=\tilde g_1$ and $\tilde g=\tilde g_1$ we get $E_N\tilde g_1=E_N\tilde g_2$. We notice that this gives us the integer $k$ of Theorem \ref{th0}.
\par\medskip
Now let $g$ be a holomorphic function on $X\cap D$ which admits a smooth extension $\tilde g$ which satisfies the assumptions of Theorem \ref{th0}. We prove that $E_N(g)(z_0)=g(z_0)$ for all $z_0\in X\cap D$.\\
For $\varepsilon>0$ small enough we construct $P_\varepsilon^{N,n}$, the Berndtsson-Andersson kernel of the domain $D_{-\varepsilon}$ which has the defining function $\rho_\varepsilon=\rho+\varepsilon$. We set $P^{N,n}_\varepsilon(\zeta ,z)=0$  for $\zeta \notin D_{-\varepsilon}$.
The kernel $P_{\varepsilon}^{N,n}(\cdot,z_0)$ converges to $P^{N,n}(\cdot,z_0)$ when $\varepsilon$ tends to $0$.

Now let 
$g_\varepsilon$ be an holomorphic extension of $g$ on $D_{-\frac\varepsilon2}$ given by Cartan's Theorem B. Fact~1 yields 
\begin{eqnarray*}
g(z_0)&=&g_{\varepsilon}(z_0)\\
&=&\int_D g_\varepsilon (\zeta )\wedge P_{\varepsilon}^{N,2}(\zeta ,z_0)\\
&=&T\left[fg_\varepsilon\wedge P_{\varepsilon}^{n,2}(\cdot,z_0)\right]\\
&=&C_1\overline\partial T\left[g_\varepsilon b(\cdot,z_0)\wedge P_{\varepsilon}^{N,1}(\cdot,z_0)\right].
\end{eqnarray*}
Then, since $P_\varepsilon^{N,1}$ is supported in $D_{-\varepsilon}$, since $\tilde g=g_\varepsilon$ on $X\cap D_{-\frac\varepsilon2}$ and since $\diffp{^{\alpha+\beta}\tilde g}{\overline{\zeta_1^*}^\alpha\partial \overline{\zeta_2^*}^\beta}=0$ on $D_{-\frac\varepsilon2}\cap X$, fact~2 gives
$$g(z_0)=C_1\overline\partial T\left[\tilde g b(\cdot,z_0)\wedge P_{\varepsilon}^{N,1}(\cdot,z_0)\right]$$
and when $\varepsilon$ goes to $0$, this yields $g(z_0)=E_N \tilde g(z_0)$ and thus $E_N g$ is an extension of $g$.
\section{Estimate of the extension operator}\mlabel{secIII}
We prove in this section that $E_N(g)$ satisfies the conclusion of Theorem \ref{th0}. For this purpose
we write $b$ in the \ko coordinates at $z_j,$ as $b(\zeta,z)=\sum_{l=1,2}
b^*_l(\zeta,z)d\zeta^*_l$ where $b^*_l(\zeta,z)=\int_0^1 \diffp f
{\zeta^*_l}(\zeta+t(z-\zeta))dt$ and 
we
prove the following estimates. We recall that for any non negative integer $j$, 
$p_j$ is the integer given by proposition \ref{propII.2.1} and 
$$I_j\hskip -1.5pt =\hskip -1.5pt\left\{i, \exists z_1^*\in \Delta_0(\kappa
|\rho(z_j)|) \text{ such that } |\alpha_i^*(z^*_1)|\leq
(2\kappa|\rho(z_j)|)^{\frac12}\right\}.$$
\begin{proposition}\mlabel{kernelestimate}
For all positive integer $j$, all $z$ in $D$ and all $\zeta$ in $\pk j$,  we have uniformly in $z,\zeta$ and $j$
\begin{eqnarray*}{\left|\frac{\prod_{i\in
I_j}\overline{\zeta^*_2-\alpha_i^*(\zeta_1)}}{f(\zeta)}b_1(\zeta,z)\right|}
&\leqs&
\sum_{0\leq\alpha+\beta\leq p_j} \delta(\zeta,z)^{\alpha+\frac\beta2}|\rho(\zeta)|^{-1-\alpha+\frac{\# I_j-\beta}2},\\
{\left|\frac{\prod_{i\in
I_j}\overline{\zeta^*_2-\alpha_i^*(\zeta_1)}}{f(\zeta)}b_2(\zeta,z)\right|}
&\leqs&
\sum_{0\leq\alpha+\beta\leq p_j} \delta(\zeta,z)^{\alpha+\frac\beta2}|\rho(\zeta)|^{-\frac12-\alpha+\frac{\# I_j-\beta}2},\\
{\left|\frac{\prod_{i\in
I_j}\overline{\zeta^*_2-\alpha_i^*(\zeta_1)}}{f(\zeta)}d_z b_1(\zeta,z)\right|}
&\leqs&
\sum_{0\leq\alpha+\beta\leq p_j} \delta(\zeta,z)^{\alpha+\frac\beta2}|\rho(\zeta)|^{-2-\alpha+\frac{\# I_j-\beta}2},\\
{\left|\frac{\prod_{i\in
I_j}\overline{\zeta^*_2-\alpha_i^*(\zeta_1)}}{f(\zeta)}d_z b_2(\zeta,z)\right|}
&\leqs&
\sum_{0\leq\alpha+\beta\leq p_j} \delta(\zeta,z)^{\alpha+\frac\beta2}|\rho(\zeta)|^{-\frac32-\alpha+\frac{\# I_j-\beta}2}.\end{eqnarray*}
\end{proposition}
\pr We prove the first inequality, the others are analogous. For $A\subset\{1,\ldots,p_j\}$ we denote by $A^c$ the complementary
of $A$ in $\{1,\ldots, p_j\}$. Proposition \ref{propII.2.1} yields:
\begin{eqnarray*}
\left|\frac{\prod_{i\in I_j}
\overline{\zeta^*_2-\alpha_i^*(\zeta^*_1)}}{f(\zeta)} \right| &\leqs&
\frac1{\prod_{i\in I_j^c} |\zeta^*_2-\alpha^*_i(\zeta_1^*)|}
\end{eqnarray*}
uniformly \wrt $\zeta$ and $j$.\\
We estimate $b_1^*$.
We have
$$\diffp{f}{\zeta^*_1}(\zeta+t(z-\zeta))=\sum_{0\leq \alpha+\beta\leq
p_j}\diffp{^{\alpha+\beta+1}f}{{\zeta^*_1}^{\alpha+1}\partial
{\zeta_2^*}^\beta}(\zeta)(z^*-\zeta^*)^{\alpha+\beta}+o(|\zeta^*-z^*|^{
p_j})$$
and
$$\left|\diffp{^{\alpha+\beta+1}f}{{\zeta^*_1}^{\alpha+1}\partial
{\zeta^*_2}^\beta}(\zeta)\right|=\left|\sum_{\over{n_1+\ldots n_{p_j}=\alpha+1}{F_1\dot{\cup}
F_2\dot{\cup}F_3=\{1,\ldots, p_j\}}}\prod_{i\in F_1}\diffp{^{n_i}
\alpha^*_i}{{\zeta^*_1}^{n_i}}(\zeta^*_1) \prod_{i\in
F_3}(\zeta^*_2-\alpha_i^*(\zeta^*_1))\right|$$
where $\dot{\cup}$ means that the union is disjoint, $F_1=\{i,\ n_i\neq 0\}$ and
$\#F_2=\beta$.\\
Since $\diffp{\alpha^*_i}{\zeta^*_1}$ is uniformly bounded and holomorphic on
$\Delta_0(2\kappa|\rho(z_j)|)$, we have 
$\left|\diffp{^{n_i}\alpha_i^*}{{\zeta^*_1}^{n_i}}\right|\leqs |\rho(z_j)|^{-n_i+1}$ on $\Delta_0(\kappa|\rho(z_j)|)$.
Moreover Proposition \ref{propII.0.1} gives $|\rho(z_j)|\eqs |\rho(\zeta)|$ for all $\zeta\in\pk
j$ so
\begin{eqnarray*}
{\left|\diffp{^{\alpha+\beta+1}f}{{\zeta^*_1}^{\alpha+1}\partial
{\zeta^*_2}^\beta}(\zeta)\right|}
&\leqs&\sum_{\over{n_1+\ldots n_{p_j}=\alpha+1}{\over{F_1\dot{\cup}
F_2\dot{\cup}F_3=\{1,\ldots, p_j\}}{\# F_2=\beta}}} |\rho(\zeta)|^{-\alpha-1+\# F_1} \prod_{i\in F_3}
|\zeta^*_2-\alpha^*_i(\zeta^*_1)|
\end{eqnarray*}
and so
\begin{eqnarray*}
{|b^*_1(\zeta,z)|}
&\leqs&\sum_{0\leq\alpha+\beta\leq p_j}\sum_{{\over{F_1\dot{\cup}
F_2\dot{\cup}F_3=\{1,\ldots, {p_j}\}}{\#F_2=\beta}}} |\rho(\zeta)|^{-1-\alpha+\#
F_1}\delta(\zeta,z)^{\alpha+\frac\beta2} \prod_{i\in F_3}
|\zeta^*_2-\alpha^*_i(\zeta^*_1)|.\\
\end{eqnarray*}
Therefore 
$\frac{\prod_{i\in I_j} \overline{\zeta^*_2-\alpha_i^*(\zeta^*_1)}}{f(\zeta)} b^*_1(\zeta,z)$ is bounded by a sum for $0\leq\alpha+\beta\leq p_j$, $F_1\dot{\cup}
F_2\dot{\cup}F_3=\{1,\ldots, {p_j}\}$, ${\#F_2=\beta}$ of 
\begin{eqnarray*}
 {S^{\alpha,\beta}_{F_1,F_2,F_3}}
&:=&\frac{\prod_{i\in F_3}|\zeta^*_2-\alpha^*_i(\zeta_1^*)|}
{\prod_{i\in I_j^c}|\zeta^*_2-\alpha^*_i(\zeta_1^*)|}|\rho(\zeta)|^{-1-\alpha+\#F_1} \delta(\zeta,z)^{\alpha+\frac\beta2}.
\end{eqnarray*}
On the one hand for $i\in I_j^c$ and $\zeta\in \pk j$ we have
$|\zeta^*_2-\alpha_i^*(\zeta^*_1)|\geqs |\rho(z_j)|^{\frac12}\eqs
|\rho(\zeta)|^{\frac12}$.
On the other hand for $i\in I_j$ and $\zeta\in \pk j$ we have
$|\zeta^*_2-\alpha_i^*(\zeta^*_1)|\leqs  |\rho(\zeta)|^{\frac12}$. Therefore,
writing $\frac{\prod_{i\in F_3}(\zeta^*_2-\alpha^*_i(\zeta_1^*))}{\prod_{i\in
I_j^c}(\zeta^*_2-\alpha^*_i(\zeta_1^*))}$ as $\frac{\prod_{i\in F_3\cap
I_j}(\zeta^*_2-\alpha^*_i(\zeta_1^*))}{\prod_{i\in I_j^c\cap
F_3^c}(\zeta^*_2-\alpha^*_i(\zeta_1^*))}\cdot \frac{\prod_{i\in F_3\cap
I_j^c}(\zeta^*_2-\alpha^*_i(\zeta_1^*))}{\prod_{i\in I_j^c\cap
F_3}(\zeta^*_2-\alpha^*_i(\zeta_1^*))}$ we get
$$S^{\alpha,\beta}_{F_1,F_2,F_3}\leqs\delta(\zeta,z)^{\alpha+\frac\beta2}
|\rho(\zeta)|^{-1-\alpha+\#F_1 +\frac{\#F_3\cap I_j - \# F_3^c\cap
I_j^c}2}.$$
The equality ${\#F_3\cap I_j - \# F_3^c\cap I_j^c}=\#I_j-\# F_3^c$ implies that  $\#F_1 +\frac{\#F_3\cap I_j - \# F_3^c\cap
I_j^c}2\geq \frac{\# I_j-\beta} 2$.\\
This gives $S^{\alpha,\beta}_{F_1,F_2,F_3}\leqs
 \delta(\zeta,z)^{\alpha+\frac\beta2}
|\rho(\zeta)|^{-1-\alpha+\frac{\#I_j-\beta}2}$
which finally yields 
\begin{eqnarray*}{\left|\frac{\prod_{i\in
I_j}\overline{\zeta^*_2-\alpha_i^*(\zeta_1)}}{f(\zeta)}b_1(\zeta,z)\right|}
&\leqs&
\sum_{0\leq\alpha+\beta\leq p_j} \delta(\zeta,z)^{\alpha+\frac\beta2}|\rho(\zeta)|^{-1-\alpha+\frac{\# I_j-\beta}2}.\qed
\end{eqnarray*}

As usually in the estimates of the Berndtsson-Andersson kernel, the main difficulty appears when we
integrate for $\zeta$ near $z$ and $z$ near $bD$. Therefore we choose $\varepsilon_0>0$ arbitrarily  small and we divide the domain of integration in two parts
: ${\cal P}_{\frac{\varepsilon_0}{2c_1}}(z)$ and $D\setminus{\cal
P}_{\frac{\varepsilon_0}{2 c_1}}(z)$ where $c_1$ is given by Proposition
\ref{propII.0.1}. In order to estimate the integral over ${\cal
P}_{\frac{\varepsilon_0}{2c_1}}(z)$, we prove the following lemma:
\begin{lemma}\mlabel{lemcov}
For all $z\in D\setminus D_{-\frac{\varepsilon_0}2}$ such that $|\rho(z)|<\frac{\varepsilon_0}{2}$, let 
$j_0$ be an integer such that $(1-c\kappa)^{-j_0}\varepsilon_0< |\rho(z)|\leq (1-c\kappa)^{-j_0-1}\varepsilon _0$ and let $z_1^{i,j},\ldots,
z_{m_{i,j}}^{i,j}$, $i\in \nn$, $j\in\zz$, be the points of the covering such
that
\begin{itemize}
 \item[-] $\rho(z^{i,j}_m)=-(1-c\kappa)^{j-j_0}\varepsilon _0$,
 \item[-] $\delta(z_m^{i,j},z)\in [i\kappa(1-c\kappa)^{j-j_0}\varepsilon _0,(i+1)\kappa(1-c\kappa)^{j-j_0}\varepsilon _0[$,
 \item[-] $\delta(z_m^{i,j},z)\leq \varepsilon_0$.
\end{itemize}
For $j\geq j_0$ let 
$i_0(j)$ be the non negative integer such that $i_0(j)\kappa(1-c\kappa)^{j-j_0}<1
\leq (1+i_0(j))\kappa(1-c\kappa)^{j-j_0}$.\\
Then 
\begin{enumerate}[(i)]
 \item \mlabel{premierpoint} ${\cal P}_{\frac{\varepsilon_0}{2c_1}}(z)\subset
\cup_{j=j_0}^{+\infty}\cup_{i=0}^{i_0(j)}\cup_{m=1}^{m_{i,j}}{\cal
P}_{\kappa|\rho(z^{i,j}_m)|}(z_m^{i,j})$,
 \item $m_{i,j}\leqs i^2$ uniformly with respect to $z_0,z,i$ and $j$.
\end{enumerate}
\end{lemma}
\pr We first prove (\ref{premierpoint}). Let $\zeta$ be a point in ${\cal P}_{\frac{\varepsilon_0}{2c_1}}(z)$.
Proposition \ref{propII.0.1} implies that $\zeta$ belongs to $D\setminus D_{-\varepsilon_0}$ so there
exists a point $\zeta_0$ of the covering such that $\zeta$ belongs to ${\cal
P}_{\kappa|\rho(\zeta_0)|}(\zeta _0)$.\\
The point $\zeta _0$ belongs to $D\setminus D_{-\varepsilon_0}$ thus there exists $j\geq j_0$ such that
$|\rho(\zeta_0)|=(1-c\kappa)^{j-j_0}\varepsilon _0$. Moreover if $\kappa $ is small enough 
\begin{eqnarray*}
 \delta(\zeta _0,z)&\leq& c_1(\delta(\zeta,\zeta_0)+\delta(\zeta,z))\\
&\leq& c_1\left(\kappa(1-c\kappa)^{j-j_0}\varepsilon _0+\frac{\varepsilon_0}{2c_1}\right)\\
&\leq&\varepsilon_0.
\end{eqnarray*}
So there exists $i\in\nn$ such that $\delta(\zeta_0,z)$ belongs to
$[i\kappa(1-c\kappa)^{j-j_0}\varepsilon _0,(i+1)\kappa(1-c\kappa)^{j-j_0}\varepsilon _0[$ and $(i+1)\kappa(1-c\kappa)^{j-j_0}\varepsilon _0\leq
\varepsilon_0$ which means that $i\leq i_0(j)$. Thus $\zeta_0$ is one the
points $z_1^{i,j},\ldots, z_{m_{i,j}}^{i,j}$ and (\ref{premierpoint}) holds.\\
In order to prove that $m_{i,j}\leqs i^2$ we introduce the set
$$E_{i,j}=\{\zeta\in D,\ \rho(\zeta)=-(1-c\kappa)^{j-j_0}\varepsilon _0\text{ and }
\delta(\zeta,z)\leq c_1\kappa(i+2)(1-c\kappa)^j|\rho(z)|\}.$$
On the one hand we have
\begin{eqnarray}
\sigma(E_{i,j})&=&\sigma\left(bD_{-(1-c\kappa)^j|\rho(z_0)|}\cap \nonumber {\cal
P}_{c_1\kappa(i+2)(1-c\kappa)^j|\rho(z)|}(z)\right)\\
&\leq& \left(c_1\kappa(i+2)(1-c\kappa)^j|\rho(z)|\right)^2\nonumber\\
&\leqs&\left(c_1\kappa(i+2)(1-c\kappa)^{j-j_0}\varepsilon _0\right)^2\mlabel{eq2}
\end{eqnarray}
On the other one hand for all $m$, all $\zeta\in {\cal
P}_{\kappa|\rho(z^{i,j}_m)|}(z_m^{i,j})$ we have:
\begin{eqnarray*}
 \delta(\zeta,z)&\leq& c_1(\delta(\zeta,z_m^{i,j})+\delta(z_m^{i,j},z))\\
&\leq& c_1(\kappa(1-c\kappa)^{j-j_0}\varepsilon _0+\kappa(i+1)(1-c\kappa)^{j-j_0}\varepsilon _0)\\
&\leq&c_1\kappa(i+2)(1-c\kappa)^{j-j_0}\varepsilon _0.
\end{eqnarray*}
This implies that ${\cal P}_{\kappa|\rho(z^{i,j}_m)|}(z_m^{i,j})\cap
bD_{-(1-c\kappa)^{j-j_0}\varepsilon _0}\subset E_{i,j}$ for all $m$ and so
\begin{eqnarray}
 \sigma(E_{i,j})&\geq& \sigma\left(\cup_{m=1}^{m_{i,j}} {\cal
P}_{\kappa|\rho(z^{i,j}_m)|}(z_m^{i,j})\cap bD_{-(1-c\kappa)^{j-j_0}\varepsilon _0} \right)\nonumber.
\end{eqnarray}
Now, the construction of a $\kappa$-covering and Proposition \ref{propII.0.1} implies that the intersection of ${\cal P}_{\frac{c\kappa}{c_1}|\rho(z^{i,j}_m)|}(z_m^{i,j})$ and ${\cal
P}_{\frac{c\kappa}{c_1}|\rho(z^{i,j}_l)|}(z_l^{i,j})$ is empty for for $l\neq m$. Therefore we have
\begin{eqnarray}
 \sigma(E_{i,j})&\geq& \sum_{m=1}^{m_{i,j}}\sigma\left({\cal P}_{\frac{c\kappa}{c_1}|\rho(z^{i,j}_m)|}(z_m^{i,j})\cap
bD_{-(1-c\kappa)^{j-j_0}\varepsilon _0} \right),\nonumber\\
&\geq& m_{i,j}(\frac{c\kappa}{c_1}(1-c\kappa)^{j-j_0}\varepsilon _0)^2.\mlabel{eq3}
\end{eqnarray}
Inequalities (\ref{eq2}) and (\ref{eq3}) together imply that $m_{i,j}\leqs
i^2$, uniformly with respect to $z$, $i$ and $j$.\qed

In order to prove the $BMO$-estimates of Theorem \ref{th0} we apply the following classical
lemma:
\begin{lemma}
 Let $h$ be a function of class $C^1$ on $D$. If there exists $C>0$ such that ${\rm d}h(\zeta)\leq
C|\rho(\zeta)|^{-1}$ then $h$ belongs to ${BMO}(D)$ and $\|h\|_{BMO(D)}\leq C$.
\end{lemma}
\noindent{\it Proof of Theorem \ref{th0} for $q=+\infty$~: } 
Let $g$ be a holomorphic function on $X\cap D$ which have a smooth extension $\tilde g$ which satisfies the assumptions (\ref{th0i}), (\ref{th0ii}) and (\ref{th0iii}) of Theorem \ref{th0}. We put 
$\gamma_\infty=\sup_{\over{\zeta\in D}{\alpha+\beta\leq k}}\left|\diffp{^{\alpha+\beta}\tilde g}{\overline{\eta_\zeta}^\alpha\partial \overline{v_\zeta}^\beta}(\zeta)\right||\rho(\zeta)|^{\alpha+\frac\beta 2}$
In order to prove Theorem \ref{th0} when $q=+\infty$, we have to prove that $E_N g$ is in $BMO(D)$ and $\|E_N g\|_{BMO(D)}\leqs \gamma_\infty$.\\
Since the Berndtsson-Andersson kernel is regular when $\zeta$ and $z$ are far from each other or when $z$ is far from $bD$, we only have to estimate the integral over ${\cal P}_{\frac{\varepsilon_0}{2c_1}}(z)$ for $z$ near $bD$ and $\varepsilon_0>0$ not depending from $z$. We keep the notation of lemma \ref{lemcov} and use the covering $\cup_{j=j_0}^{+\infty}\cup_{i=0}^{i_0(j)}\cup_{m=1}^{m_{i,j}}{\cal
P}_{\kappa|\rho(z^{i,j}_m)|}(z_m^{i,j})$ of ${\cal P}_{\frac{\varepsilon_0}{2c_1}}(z)$ given by lemma \ref{lemcov}.
We denote by $p_m^{i,j}$ the number of sheets given by proposition \ref{propII.2.1} for $z_m^{i,j}$, $I_{m}^{i,j}$ is the set
$I_m^{i,j}\hskip -1.5pt =\hskip -1.5pt\left\{k, \exists z_1^*\in \Delta_0(\kappa|\rho(z_m^{i,j})|) \text{ such that } |\alpha_k^*(z^*_1)|\leq
(2\kappa|\rho(z^{i,j}_m)|)^{\frac12}\right\}$ and $q_m^{i,j}$ denotes its cardinal.\\
From Proposition \ref{estiBA} and \ref{kernelestimate} we get for all $\zeta\in {\cal
P}_{\kappa|\rho(z^{i,j}_m)|}(z_m^{i,j})$
\begin{eqnarray*}
 \lefteqn{\left|{d}_z \left(\frac{\prod_{i\in
I_m^{i,j}}\overline{\zeta^*_2-\alpha_i^*(\zeta_1)}}{f(\zeta)}b(\zeta,z)\wedge \overline \partial\diffp{^{q^{i,j}_m}}{\overline{\zeta^*_2}^{q^{i,j}_m}} \left(\tilde g(\zeta)P^{N,n}(\zeta,z)\right) \right)\right|}\\
&\leqs& \gamma_\infty\sum_{0\leq \alpha+\beta\leq p_m^{i,j}} 
\left(\frac{\delta(\zeta,z)}{|\rho(\zeta)|}\right)^{\alpha+\frac\beta2} \frac{|\rho(\zeta)|^{N}}{(|\rho(\zeta)|+|\rho(z)|+\delta(z,\zeta))^{N+4}}
\\
&\leqs&\gamma_\infty 
\frac{|\rho(\zeta)|^{N'}}{(|\rho(\zeta)|+|\rho(z)|+\delta(z,\zeta))^{N'+4}}.
\end{eqnarray*}
where $N'=N-\max_{i,j} p_{i,j}$.\\
We have for all $\zeta\in {\cal P}_{\kappa|\rho(z^{i,j}_m)|}(z_m^{i,j})$, $|\rho(\zeta)|\geq \frac12 |\rho(z^{i,j}_m)|$ and thus:
\begin{eqnarray*}
|\rho(\zeta)|+\delta(\zeta,z)
&\geq& \frac12|\rho(z^{i,j}_m)|+\frac1{c_1} \delta(z,z_m^{i,j})-\delta(z^{i,j}_m,\zeta)\\
&\geq& |\rho(z^{i,j}_m)|(\frac12-\kappa)+\frac1{c_1}\delta(z,z_m^{i,j})\\
&\geqs& |\rho(z_{m}^{i,j})|+\delta(z,z_m^{i,j}).
\end{eqnarray*}
Therefore
 \begin{eqnarray*}
 \lefteqn{\left|{d}_z \left( \frac{\prod_{i\in
I_m^{i,j}}\overline{\zeta^*_2-\alpha_i^*(\zeta_1)}}{f(\zeta)}b(\zeta,z)\wedge \overline \partial\diffp{^{q^{i,j}_m}}{\overline{\zeta^*_2}^{q^{i,j}_m}} \left(\tilde g(\zeta)P^{N,n}(\zeta,z)\right)\right) \right|}\\
&&\hskip 90pt\leqs \gamma_\infty
\frac{|\rho(z^{i,j}_m)|^{N'}}{(|\rho(z)|+|\rho(z^{i,j}_m)|+\delta(z,z^{i,j}_m))^{N'+4}}.
\end{eqnarray*}
Now, integrating over ${\cal P}_{\kappa|\rho(z^{i,j}_m)|}(z^{i,j}_m)$ and summing over $m$, $i$ and $j$ we have to prove that
the sum
$$\sum_{j=j_0}^\infty\sum_{i=0}^{i_0(j)}\sum_{m=1}^{m_{i,j}}
\frac{|\rho(z_{m}^{i,j})|^{N'}}{\left((i+1)|\rho(z^{i,j}_m)|+|\rho(z)|\right)^{N'+1}}
$$
is uniformly bounded by $\frac1{|\rho(z)|}$. We have:
\begin{eqnarray*}
\lefteqn{\sum_{j=j_0}^\infty\sum_{i=0}^{i_0(j)}\sum_{m=1}^{m_{i,j}}
\frac{|\rho(z_{m}^{i,j})|^{N'}}{\left((i+1)|\rho(z^{i,j}_m)|+|\rho(z)|\right)^{N'+1}}}\\
&\leq&
 \sum_{j=j_0}^\infty\sum_{i=0}^{i_0(j)}\sum_{m=1}^{m_{i,j}}
\left(\frac{(1-c\kappa)^j}{(i+1)(1-c\kappa)^j+1}\right)^{N'} \cdot\frac{1}{((i+1)(1-c\kappa)^j+1)|\rho(z)|}\\
&\leq&\frac{1}{|\rho(z)|}\left(\sum_{j=0}^\infty \sum_{i=0}^\infty
\frac{(1-c\kappa)^j}{(i+1)^{N'-3}} +
\sum_{j=j_0}^{-1} \sum_{i=0}^\infty \frac1{(i+1)^{N'-2}(1-c\kappa)^{j}}\right)\\
&\leqs& \frac1{|\rho(z)|}.
\end{eqnarray*}
So $E_N(g)$ belongs to ${BMO}(D)$ and $\|E_N(g)\|_{BMO(D)} \leqs \sup_{\over{\zeta\in D}{\alpha+\beta\leq k}}\left|\diffp{^{\alpha+\beta}\tilde g}{\overline{\eta_\zeta}^\alpha\partial \overline{v_\zeta}^\beta}(\zeta)\right||\rho(\zeta)|^{\alpha+\frac\beta 2}$.\qed\\
The $L^q$-estimates of Theorem \ref{th0} are left to be shown. For $q\in (1,+\infty)$ we will apply the following lemma (see
\cite{Pol}):
\begin{lemma}\mlabel{Pol}
 Suppose the kernel $k(\zeta,z)$ is defined on $D\times D$ and the operator $K$
is defined by $Kf(z)=\int_{\zeta\in D}k(\zeta,z)f(\zeta)d\lambda(\zeta)$. If for
every $\varepsilon\in ]0,1[$ there exists a constant $c_\varepsilon$ such that
\begin{eqnarray*}
 \int_{\zeta\in D} |\rho(\zeta)|^{-\varepsilon}|k(\zeta,z)|d\lambda(\zeta)&\leq&
c_\varepsilon |\rho(z)|^{-\varepsilon},\quad \forall z\in D,\\
\int_{z\in D} |\rho(z)|^{-\varepsilon}|k(\zeta,z)|d\lambda(z)&\leq& c_\varepsilon
|\rho(\zeta)|^{-\varepsilon},\quad \forall \zeta\in D
\end{eqnarray*}
Then for all $q\in ]1,+\infty[$, there exists $c_q>0$ such that
$\|Kf\|_{L^q(D)}\leq \|f\|_{L^q(D)}$.
\end{lemma}
{\it Proof of Theorem \ref{th0} for $q\in(1,+\infty)$~:}
Applying Lemma \ref{Pol} and Propositions \ref{estiBA} and  \ref{kernelestimate}, it
suffices to prove that for all $\varepsilon\in (0,1)$ there exists
$c_\varepsilon>0$ such that
\begin{eqnarray}
{\int_{\zeta\in D}
\frac {|\rho(\zeta)|^{N'-\varepsilon}}
{\left(|\rho(\zeta)|+|\rho(z)|+\delta(\zeta ,z)\right)^{N'+3}}
d\lambda(\zeta)}&\leq &c_{\varepsilon}|\rho(z)|^{-\varepsilon},\ \forall z\in D,\mlabel{eq5}\\
{\int_{z\in D}
\frac {|\rho(\zeta)|^{N'}|\rho(z)|^{-\varepsilon}}
{\left(|\rho(\zeta)|+|\rho(z)|+\delta(\zeta ,z)\right)^{N'+3}}
d\lambda(z)}&\leq &c_{\varepsilon}|\rho(\zeta)|^{-\varepsilon},\ \forall \zeta \in D,\mlabel{eq6}
\end{eqnarray}
The inequality (\ref{eq5}) can be shown as in the proof of Theorem \ref{th0} for $q=\infty$.\\
In order to prove that the inequality (\ref{eq6}) holds true 
we cover $D$ with the \ko balls ${\cal P}_{\kappa|\rho(\zeta)|}(\zeta)$ and $\left({\cal P}_{2^{j+1}\kappa|\rho(\zeta)|}(\zeta)\setminus {\cal
P}_{2^j\kappa|\rho(\zeta)|}(\zeta)\right)$, $j\in\nn$.\\
For $z\in{\cal P}_{\kappa|\rho(\zeta)|}(\zeta )$, $|\rho(z)|\eqs |\rho(\zeta)|$ and thus
\begin{eqnarray}
\int_{z\in{\cal P}_{\kappa|\rho(\zeta)|}(\zeta )} \frac {|\rho(\zeta)|^{N'}|\rho(z)|^{-\varepsilon}}
{\left(|\rho(\zeta)|+|\rho(z)|+\delta(\zeta ,z)\right)^{N'+3}}d\lambda(z)
&\leqs&|\rho(\zeta)|^{-\varepsilon}.\mlabel{eq20}
\end{eqnarray}
When we integrate on ${\cal P}_{2^{j+1}\kappa|\rho(\zeta)|}(\zeta)\setminus {\cal
P}_{2^j\kappa|\rho(\zeta)|}(\zeta)$ we get
\begin{eqnarray}
\nonumber\lefteqn{ \int_{{\cal P}_{2^{j+1}\kappa |\rho(\zeta)|}(\zeta)\setminus {\cal
P}_{2^j\kappa |\rho(\zeta)|}(\zeta)} 
 \frac {|\rho(\zeta)|^{N'}|\rho(z)|^{-\varepsilon}}
{\left(|\rho(\zeta)|+|\rho(z)|+\delta(\zeta ,z)\right)^{N'+3}}d\lambda(z)}\\
\nonumber&\hskip 100pt &\leqs
\int_{\over{|x_1|,|y_1|\leq 2^{j+1}\kappa|\rho(\zeta)|}{|x_2|,|y_2|\leq \sqrt{2^{j+1}\kappa|\rho(\zeta)|}}}
 \frac{|\rho(\zeta)|^{N'}x_1^{-\varepsilon}}
{\left(|\rho(\zeta)|+2^j\kappa|\rho(\zeta)|\right)^{N'+3}}d\lambda(z)\\
\nonumber&\hskip 100pt &\leqs
(2^{j+1}\kappa|\rho(\zeta)|)^{-\varepsilon+3}
\frac{|\rho(\zeta)|^{N'}}
{\left(|\rho(\zeta)|+2^j\kappa|\rho(\zeta)|\right)^{N'+3}}\\
&\hskip 100pt &\leqs |\rho(\zeta)|^{-\varepsilon} 2^{-j(N'+\varepsilon)}\mlabel{eq21}
\end{eqnarray}
Summing  (\ref{eq20}) and (\ref{eq21}) for all non-negative integer $j$ we prove inequality (\ref{eq20}).  Theorem \ref{th0} is therefore proved
for $q\in(1,+\infty)$.\qed\\
{\it Proof of Theorem \ref{th0} for $q=1$~:} We prove directly that $E_Ng$ belongs to $L^1(D)$. 
Propositions \ref{estiBA} and  \ref{kernelestimate} yield
\begin{eqnarray*}
{\int_D|E_N g(z)|d\lambda(z)}&\leqs&\sum_{j=0}^\infty\sum_{0\leq\alpha+\beta\leq q_j+1} \int_{\pk j} {|\rho(z_j)|^{\alpha+\frac\beta2}} \left| \diffp{^{\alpha+\beta}\tilde g} {\overline{\zeta^*_1}^\alpha\partial \overline{\zeta^*_2}^\beta}(\zeta)\right|\\
&&\hskip 50pt \cdot \left(\int_D
\frac{|\rho(\zeta)|^{N'}}{\left(|\rho(\zeta)|+|\rho(z)|+\delta(\zeta,z)\right)^{N'+3}}
{d\lambda(z)}\right)d\lambda(\zeta).
\end{eqnarray*}
As for the proof of (\ref{eq6}) we cover $D$ using \ko corona and get
\begin{eqnarray*}
{\int_D|Eg(z)|d\lambda(z)}
&\leqs&\sum_{j=0}^\infty\sum_{0\leq \alpha+\beta\leq q_j+1} \int_{\pk j} {|\rho(z_j)|^{\alpha+\frac\beta2}} \left| \diffp{^{\alpha+\beta}\tilde g} {\overline{\zeta^*_1}^\alpha\partial \overline{\zeta^*_2}^\beta}(\zeta)\right|d\lambda(\zeta)\\
&\leqs& \sum_{0\leq\alpha +\beta \leq k}\left\| \zeta \mapsto \diffp{^{\alpha +\beta }\tilde g}{\overline\eta_\zeta ^\alpha \partial\overline v_\zeta ^\beta }(\zeta )\rho (\zeta )^{\alpha +\frac\beta 2}\right\|_{L^1(D)}.
\end{eqnarray*}
\qed
\section{Smooth extension and divided differences}\mlabel{section5}
In this section we give necessary conditions in $\cc^n$ that a function $g$ holomorphic on $X\cap D$ has to satisfy in order to have a $L^q$-holomorphic extension on $D$, $q\in [1,+\infty]$. We also prove that these conditions are sufficient in $\cc^2$  for $g$ to have a $L^q$-holomorphic extension on $D$ when $q$ belongs to $[1,+\infty)$ or a $BMO$-holomorphic extension when $q=+\infty$.
\subsection{$L^\infty$-$BMO$ extension}
We first prove the following lemma for functions defined on $X\cap D$ which have holomorphic extension on $D$. We use the notations defined in the introduction.
\begin{lemma}\mlabel{lemma0}
If $g$ defined on $X\cap D$ has a holomorphic extension $G$ on $D$ then uniformly with respect to $g$,  $G$,  $z\in D$,  $v$ unit vector of $\cc^n$ and  positive integer $k$ such that $k\leq \# \Lambda(z,v)$ :
$$\sup_{
\genfrac{}{}{0pt}{}{\lambda_1,\ldots,\lambda_k\in\Lambda_{z,v}}{\lambda_i\neq\lambda_j\text{ for } i\neq j}}|g_{z,v}[\lambda_1,\ldots,\lambda_k]|  \tau(z,v,|\rho(z)|)^{k-1} \leqs \sup_{b\Delta_{z,v} \left(4\kappa\tau(z,v,|\rho(z)|)\right)}|G|.$$
\end{lemma}
\pr For $\lambda_1,\ldots,\lambda_k\in \Lambda_{\zeta,v}$ pairwise distincts, we have by Cauchy's formula
$$g_{z,v}[\lambda_1,\ldots,\lambda_k]=\frac1{2i\pi}\int_{|\lambda|=4 \tau(z,v,|\rho(z)|)} \frac{G(z+\lambda v)}{\prod_{l=1}^k(\lambda-\lambda_i)}d\lambda.$$
since for all $\lambda_i$ we have $|\lambda_i|\leq 3\tau(z,v,|\rho(z)|)$, we get
$$|g_{z,v}[\lambda_1,\ldots,\lambda_k]\leqs \left(\frac{1}{\tauzv}\right)^{k-1}\sup_{b\Delta_{z,v}\left(4\kappa \tauzv\right)}|G|.$$
\qed

\noindent{\it Proof of Theorem \ref{th1}~:} Lemma \ref{lemma0} implies directly that $c_\infty(g)\leqs \|G\|_{L^\infty(D)}$.\qed
\par\medskip
Now we prove that an even weaker assumption than $c_\infty(g)<\infty$  is actually sufficient in $\cc^2$ for $g$ to have a smooth extension which satisfies the hypothesis of Theorem \ref{th0} for $q=\infty$ and thus for $g$ to have a holomorphic $BMO$ extension on $D$. We define for $\kappa$ and $\varepsilon_0$ positive real number
$$c^{(\infty)}_{\kappa,\varepsilon_0}(g)=\sup|g_{\zeta+z^*_1\eta_\zeta,v_\zeta}[\lambda_1,\ldots,\lambda_k]|\tau(\zeta,v_\zeta,|\rho(\zeta)|)^{k-1}$$
where the supremum is taken over $\zeta\in D\setminus D_{-\varepsilon_0}$, $z_1^*\in \cc$ such that $|z^*_1|\leq \kappa|\rho(\zeta)|$, $\lambda_1,\ldots, \lambda_k\in \Lambda_{\zeta+z^*_1\eta_\zeta,v_\zeta}$ pairwise distinct. Of course, $c^{(\infty)}_{\kappa,\varepsilon_0}(g)\leq c_\infty(g)$ and it may be simpler to check that $c^{(\infty)}_{\kappa,\varepsilon_0}(g)$ is finite than to check that $c_\infty(g)$ is finite. Moreover, as told by the following lemma, when $c^{(\infty)}_{\kappa,\varepsilon_0}(g)$ is finite, $g$ admits a smooth extension which satisfies the assumptions of Theorem \ref{th0}. 

\begin{lemma}\mlabel{lemma2}
In $\cc^2$, let $g\in \oo(X\cap D)$ be such that $c^{(\infty)}_{\kappa,\varepsilon_0} (g)<\infty$. Then there exist a neighborhood $\cal U$ of $bD$ and  $\tilde g\in C^\infty(D\cap {\cal U})$ such that
\begin{enumerate}[(i)]
 \item for all non negative integer $N$, $|\rho|^{N+1} \tilde g$ vanishes to order $N$ on $bD$,
 \item for all $\alpha$ and $\beta$ non negative integer, $\left|\diffp{^{\alpha+\beta}\tilde g}{\overline{\eta_\zeta}^\alpha\partial \overline{v_\zeta}^\beta}\right||\rho|^{\alpha+\frac\beta 2}$ is bounded up to a uniform multiplicative constant on $D\cap{\cal U}$  by $c^{(\infty)}_{\kappa,\varepsilon_0}(g)$ ,
 \item for all $\alpha$ and $\beta$ non negative integer, $\diffp{^{\alpha+\beta}\tilde g}{\overline{\eta_\zeta}^\alpha\partial \overline{v_\zeta}^\beta}=0$ on $X\cap D\cap{\cal U}$.
\end{enumerate}
\end{lemma}
\pr For $\varepsilon_0>0$, we cover $D\setminus D_{-\varepsilon_0}$ with a $\kappa $-covering $\left({\cal P}_{\kappa|\rho(z_j)|}(z_j)\right)_{j\in\nn}$ constructed in subsection \ref{maxcover}. For a fixed nonnegative integer $j$, we set $w_1^*=\eta_{z_j}$ and $w^*_2=v_{z_j}$. Let $\alpha_1,\ldots,\alpha_{p_j}$ be the parametrization given by proposition \ref{propII.2.1}, $I_j=\{i,\ \exists z^*_1\in\cc\text{ with } |z_1^*|<\kappa |\rho(z_j)| \text{ and } |\alpha_i(z^*_1)| \leq 2\kappa|\rho(z_j)|\}$, $q_j=\# I_j$.\\
If $I_j=\emptyset$ we put $\tilde g_j=0$ on ${\cal P}_{\kappa|\rho(z_j)|}(z_j)$.\\
Otherwise, without restriction we assume that $I_j=\{1,\ldots, q_j\}$ and for $z=z_j+z^*_1w_1^*+z^*_2w^*_2\in {\cal P}_{2\kappa|\rho(z_j)|}(z_j)$ we put

$$\tilde g_j(z)
=\sum_{k=1}^{q_j} g_{z_j+z^*_1w^*_1,w_2^*} [\alpha_1(z^*_1),\ldots, \alpha_k (z^*_1)]\prod_{l=1}^{k-1}(\zeta^*_2- \alpha_{l}( z^*_1)).$$

Proposition \ref{propII.2.1} implies  for all $z^*_1\in \Delta_0(2\kappa |\rho (z_j)|)$ that $\alpha_j(z^*_1)$ belongs to $\Lambda_{z_j+z^*_1w^*_1,w^*_2}$ thus $\tilde g_j$ is well defined on ${\cal P}_{2\kappa|\rho(z_j)|}(z_j)$.\\
The function $\zeta\mapsto \tilde g_j(z_j+z^*_1w^*_1 +\zeta w^*_2)$ is the polynomial which interpolates $\zeta\mapsto g(z_j+z^*_1w^*_1 +\zeta w^*_2)$ at the points $\alpha_1(z^*_1),\ldots, \alpha_{q_j}(z^*_1)$ and thus $\tilde g_j$ is a holomorphic extension of $g$ on ${\cal P}_{\kappa|\rho(z_j)|}(z_j)$.\\
For all $z=z_j+z^*_1w_1^*+z^*_2w^*_2\in {\cal P}_{2\kappa|\rho(z_j)|}(z_j)$, we have $$|z^*_2-\alpha_l(z^*_1)| \leq \tau(z_j,w^*_2,2\kappa|\rho(z_j)|)\leqs \tau(z,w^*_2,2\kappa|\rho(z)|)$$  thus $|\tilde g_j(z)|\leqs c_\infty(g)$ on ${\cal P}_{2\kappa|\rho(z_j)|}(z_j)$ and $|\rho(z_j)|^{\alpha+\frac\beta2} \left|\diffp{^{\alpha+\beta} \tilde g_j}{{w^*_1}^\alpha\partial{w^*_2}^\beta}(z)\right|\leqs c_\infty(g)$ on ${\cal P}_{\kappa|\rho(z_j)|}(z_j)$.
Now we glue together all the $\tilde g_j$ using a suitable partition of unity and get our extension on $D\setminus D_{-\varepsilon_0}$. Let $(\chi_j)_{j\in\nn}$ be a partition of unity subordinated to
$\left(\pk j\right)_{j\in\nn}$ such that 
for all $j$ and all
$\zeta\in\pk j$, we have
$\left|\diffp{^{\alpha+\overline\alpha+\beta+\overline\beta} \chi_j}
{{w^{*}_1}^\alpha\partial {w^{*}_2}^\beta\partial
\overline{w^{*}_1 }^{\overline\alpha}\partial
\overline{w^{*}_2 }^{\overline\beta}}
(\zeta)\right|\leqs \frac{1}{|\rho(z_j)|^{\alpha+\overline\alpha+\frac{\beta+\overline\beta}{2}}}$, uniformly with respect to $z_j$ and $\zeta$.\\
We set $\tilde g_{\varepsilon_0}=\sum_j\chi_j\tilde g_j$. By construction for all $N\in\nn$, $\rho^{N+1} \tilde g_{\varepsilon_0}$ is of class $C^{N}$ on $\overline D\setminus D_{-\varepsilon_0}$ and vanishes to order $N$ on $bD$. Moreover, since for all $j$ the function $\tilde g_j$ is holomorphic, 
$\diffp{^{\alpha+\beta}\tilde g_{\varepsilon_0}}{\overline z^\alpha_1\partial \overline z^\beta_2}=0$ on $X\cap ( D\setminus D_{-\varepsilon})$ and, by our choice of $\chi_j$, 
$\left|\diffp{^{\alpha+\beta}\tilde g_{\varepsilon_0}}{\overline{\eta_\zeta}^\alpha\partial \overline{v_\zeta}^\beta}(\zeta )\right|\leqs |\rho(\zeta )|^{-\left(\alpha+\frac\beta 2\right)}$ for all $\zeta \in D\setminus D_{-\varepsilon _0}$.\qed\\
As a direct corollary of Lemma \ref{lemma2}, we have
\begin{corollary}\mlabel{th2}
In $\cc^2$, let $g\in \oo(X\cap D)$ be such that $c_\infty(g)<\infty$. Then there exist a neighborhood $\cal U$ of $bD$ and  $\tilde g\in C^\infty(D\cap {\cal U})$ such that
\begin{enumerate}[(i)]
 \item for all non negative integer $N$, $|\rho|^{N+1} \tilde g$ vanishes to order $N$ on $bD$,
 \item for all $\alpha$ and $\beta$ non negative integer, $\left|\diffp{^{\alpha+\beta}\tilde g}{\overline{\eta_\zeta}^\alpha\partial \overline{v_\zeta}^\beta}\right||\rho|^{\alpha+\frac\beta 2}$ is bounded up to a uniform multiplicative constant on $D\cap{\cal U}$  by $c_\infty(g)$ ,
 \item for all $\alpha$ and $\beta$ non negative integer, $\diffp{^{\alpha+\beta}\tilde g}{\overline{\eta_\zeta}^\alpha\partial \overline{v_\zeta}^\beta}=0$ on $X\cap D\cap{\cal U}$.
\end{enumerate}
\end{corollary}
Theorem \ref{th3} is now a corollary of Theorem \ref{th0} and Corollary \ref{th2} :\\
{\it Proof of Theorem \ref{th3}~:} We use Corollary \ref{th2} to get an extension $\tilde g$ of $g$ which satisfies the hypothesis of Theorem \ref{th0} on ${\cal U}\cap D$. Cartan's Theorem B gives us a bounded holomorphic extension on $D\setminus {\cal U}$. Gluing these two extensions together, we get a smooth extension of $g$ which satisfies the hypothesis of Theorem \ref{th0} in the whole domain $D$ and thus, Theorem \ref{th0} ensure the existence of a $BMO$ holomorphic extension of $g$.\qed
\subsection{$L^q(D)$-extension}
The case of $L^q$-extensions is a bit harder to handle because it is not a punctual estimate but an average estimate. Therefore the assumption under which a function $g$ holomorphic on $X\cap D$ admits a $L^q$-holomorphic extension on $D$ uses a $\kappa $-covering $\left({\cal P}_{\kappa |\rho(z_j)|}(z_j)\right)_{j\in\nn}$ in addition to the divided differences.\\
By transversality of $X$ and $bD$, for all $j$ there exists $w_j$ in the complex tangent plane to $bD_{\rho(z_j)}$ such that $\pi_j$, the orthogonal projection on the hyperplane orthogonal to $w_j$ passing through $z_j$, is a $p_j$ sheeted covering of $X$. We denote by $w_1^*,\ldots, w^*_n$ an orthonormal basis of $\cc^n$ such that $w_1^*=\eta_{z_j}$ and $w_n^*=w_j$ and we set ${\cal P}'_{\varepsilon}(z_j)=\{z'=z_j+z^*_1w^*_1+\ldots+z^*_{n-1} w^*_{n-1},\ |z^*_1|< \varepsilon \text{ and } |z_k^*|<\varepsilon^{\frac12},\ k=2,\ldots, n-1\}$. 
We put
\begin{eqnarray*}
c^{(q)}_{\kappa,{(z_j)_{j\in\nn}}}(g)\hskip-1pt &=&\hskip-1pt\sum_{j=0}^\infty \int_{z'\in{\cal P}'_{2\kappa |\rho(z_j)|}(z_j)}
\sum_{\over{\lambda_1,\ldots,\lambda_k\in\Lambda_{z',w_n^*}}{\lambda_i\neq\lambda_l\text{ for }i\neq l}}\hskip - 3pt
|\rho(z_j)|^{q\frac{k-1}2+1} \left|g_{z',w_n^*}[\lambda_1,\ldots,\lambda_k]\right| dV_{n-1}(z')
\end{eqnarray*}
where $dV_{n-1}$ is the Lebesgue measure in $\cc^{n-1}$.
\begin{theorem}\mlabel{th4}
In $\cc^n$, $n\geq 2$, let $\left({\cal P}_{\kappa |\rho(z_j)|}(z_j)\right)_{j\in\nn}$ be a $\kappa $-covering of $D\cap X$. If $g\in\oo(X\cap D)$ has a holomorphic extension $G\in L^q(D)$ then $c^{(q)}_{\kappa,(z_j)_{j\in\nn}}(g)\leqs \|G\|^q_{L^q(D)}$ uniformly with respect to $g$, $G$ and the covering $\left({\cal P}_{\kappa |\rho(z_j)|}(z_j)\right)_{j\in\nn}$.
\end{theorem}
\pr For all $j\in\nn$ all $z'\in{\cal P}_{\kappa|\rho(z_j)|}(z_j)$, all $r\in\rr$ such that $\frac72 \kappa |\rho(z_j)|^{\frac12}\leq r\leq 4\kappa |\rho(z_j)|^{\frac12}$, all $\lambda_1,\ldots, \lambda_k\in\Lambda_{z',w_n^*}$ pairwise distinct we have by Cauchy's formula
$$g_{z',w_j}[\lambda_1,\ldots,\lambda_k]=\frac1{2i\pi}\int_{|\lambda|=r} \frac{G(z'+\lambda w_j)}{\prod_{l=1}^k(\lambda-\lambda_i)}d\lambda.$$
After integration for $r\in[7/2\kappa |\rho (z_j)|,4\kappa |\rho (z_j)|]$, 
Jensen's inequality yields
$$\left|g_{z',w_j}[\lambda_1,\ldots,\lambda_k]\right|^q\leqs |\rho(z_j)|^{\frac{1-k}2 q-1} \int_{|\lambda|\leq (4\kappa |\rho(z_j)|)^\frac12} |G(z'+\lambda w_j)|^q dV_1(\lambda)$$
and thus
\begin{eqnarray*}
{\int_{z'\in {\cal P}'_{\kappa|\rho(z_j)|}(z_j)} \left|g_{z',w_j}[\lambda_1,\ldots,\lambda_k]\right|^q |\rho(z_j)|^{\frac{k-1}2 q-1}dV_{n-1}}
& \leqs& \int_{z\in {\cal P}_{4\kappa|\rho(z_j)|} (z_j)} |G(z)|^q dV_n(\lambda). 
\end{eqnarray*}
Since $\left({\cal P}_{\kappa |\rho(z_j)|}(z_j)\right)_{j\in\nn}$ is a $\kappa $-covering, we deduce from this inequality that 
$c^{(q)}_{\kappa,(z_j)_{j\in\nn}}(g)\leqs \|G\|^q_{L^q(D)}$.\qed\\[5pt]
Now we come back in $\cc^2$ and prove that the condition $c^{(q)}_{\kappa,(z_j)_{j\in\nn}}(g)<\infty$ is indeed sufficient for $g$ to have a $L^q$ extension.
\begin{theorem}
 \mlabel{th5}
In $\cc^n2,$ let $\left({\cal P}_{\kappa |\rho(z_j)|}(z_j)\right)_{j\in\nn}$ be a $\kappa $-covering of $D\cap X$. If the function holomorphic on $X\cap D$ satisfies is such that $c^{(q)}_{\kappa,(z_j)_{j\in\nn}}(g)<\infty$, then there exist a neighborhood $\cal U$ of $bD$ and a smooth extension  $\tilde g\in C^\infty(D\cap {\cal U})$ of $g$ such that
\begin{enumerate}[(i)]
 \item for all $N\in\nn$ such that $|\rho|^{N+4} \tilde g$ vanishes to order $N$ on $bD$,
 \item for all non negative integers $\alpha $ and $\beta $ the function $\zeta\mapsto\left|\diffp{^{\alpha+\beta}\tilde g}{\overline{\eta_\zeta}^\alpha\partial \overline{v_\zeta}^\beta}(\zeta )\right||\rho(\zeta )|^{\alpha+\frac\beta 2}$ has a $L^q$ norm on $D\cap {\cal U}$  bounded  by $c^{(q)}_{\kappa,(z_j)_{j\in\nn}}(g)$ up to a uniform multiplicative constant,
 \item for all non negative integer $\alpha$ and $\beta$, $\diffp{^{\alpha+\beta}\tilde g}{\overline{\eta_\zeta}^\alpha\partial \overline{v_\zeta}^\beta}=0$ on $X\cap D\cap{\cal U}$.
\end{enumerate}
\end{theorem}
\pr We proceed as in the proof of Theorem \ref{th0}. Let $\varepsilon_0$ be a positive real number. On $D\setminus D_{-\varepsilon_0}$ we define, for any non negative integer $j$,  $\chi_j$ and $\tilde g_j$ and $\tilde g_{\varepsilon _0}$  as in the proof of Lemma \ref{lemma2} and prove that it satisfies the wanted estimates. As in the proof of Lemma \ref{lemma2}, $\rho^{N+4}\tilde g_{\varepsilon_0}$ vanishes at order $N$ on $bD$ and $\diffp{^{\alpha+\beta} \tilde g_{\varepsilon_0}}{\overline{z_1}^\alpha\partial \overline {z_2}^\beta}=0$ on $X\cap D$. Moreover we have for $z\in{\cal P}_{\kappa|\rho(z_j)|}(z_j)$
\begin{eqnarray*}
 \left|
\tilde g_j(z) \diffp{^{\alpha+\beta} \chi_j}{\overline{\eta_z}^\alpha\partial \overline{v_z}^\beta} (z)\right|
&\leqs& |\rho(z_j)|^{-\alpha-\frac\beta2}\left|\tilde g_j(z)\right|\\
&\leqs& |\rho(z_j)|^{-\alpha-\frac\beta2}\sum_{k=1}^{q_j} \left|g_{z_j,v_{z_j}} [\alpha_1(z^*_1),\ldots, \alpha_k(z^*_1)]\right| |\rho(z_j)|^{\frac{k-1}2}\\
&\leqs& |\rho(z)|^{-\alpha-\frac\beta2}\sum_{k=1}^{q_j} \left|g_{z_j,v_{z_j}} [\alpha_1(z^*_1),\ldots, \alpha_k(z^*_1)]\right| |\rho(z)|^{\frac{k-1}2}
\end{eqnarray*}
and thus $z\mapsto |\rho(z)|^{\alpha+\frac\beta2} \diffp{^{\alpha+\beta} \tilde g_{\varepsilon_0}}{\overline{\eta_z}^\alpha\partial \overline {v_z}^\beta}(z)$ is in $L^q(D)$ for all $\alpha$ and $\beta$.\qed

As a corollary of Theorem \ref{th0} and Theorem \ref{th5} we get
\begin{theorem}\mlabel{th6}
In $\cc^2$, if the function $g$ holomorphic in $X\cap D$ is such that $c^{(q)}_{\kappa,(z_j)_{j\in\nn}} (g)<\infty$, then $g$ has a holomorphic extension $G$ which belongs to $L^q(D)$.
\end{theorem}
\pr Theorem \ref{th5} and Cartan's Theorem B gives a smooth extension to which we can apply Theorem \ref{th0} and get a holomorphic extension in $L^q(D)$.\qed
\subsection{Extension and weak holomorphy}
One may notice that each time, the smooth extension near the boundary is controlled only by the values of $g$ on $X\cap D$. Moreover we have never used the strong holomorphy of $g$ excepted when we involved Cartan's Theorem B in order to get a bounded extension far from the boundary. Actually, we can use only weak holomorphy and get a smooth extension and then apply theorem \ref{th0} in order to get a holomorphic extension with $BMO$ or $L^q$ norm controlled only by the values of $g$ on $X\cap D$.
Let us first recall the definition of weak holomorphy we shall use
\begin{definition}
 Let ${\cal U}$ be an open set of $\cc^n$. A function $g$ defined on $X$ is said to be weakly holomorphic on $X\cap {\cal U}$ if it is locally bounded on $X\cap{\cal U}$ and holomorphic on the regular set of $X\cap{\cal U}$.
\end{definition}
The following theorem is a direct corollary of Lemma \ref{lemma0} 
\begin{theorem}\mlabel{th7}
In $\cc^n$, for $q\in[1,+\infty)$, if the function $g$, defined on $X\cap D$, has a holomorphic extension $G\in L^q(D)$ then 
$$\sup\left|g_{z,v}[\lambda_1,\ldots,\lambda_k]\right| \tau(z,v,|\rho(z)|)^{k-1}
\left({\rm Vol}\: {\cal P}_{\kappa |\rho(z)|} (z)\right)^{\frac12} \leq \|G\|_{L^q\left({\cal P}_{\kappa |\rho(z)|} (z)\right)}$$
where the supremum is taken over all $z\in D,$ all unit vector $v$ in $\cc^n$, all positive integer $k$ such that $k\leq \#\Lambda_{z,v}$ and all $\lambda_1,\ldots,\lambda_k\in\Lambda_{z,v}$ pairwise distinct.
\end{theorem}
When $z$ is far from $bD$, Theorem \ref{th7} essentially says that the divided differences have to be bounded even in the case of $L^q$ extensions, $q<\infty$.
This is sufficient when $n=2$ to construct a smooth bounded extension in $D_{-\varepsilon}$ for $\varepsilon>0$.
\begin{lemma}\mlabel{lemma1}
 For $X$ and $D$ in $\cc^2$, let $\varepsilon$ be a positive real number. Let $g$ be a weakly holomorphic function on $X\cap D$  such that  
$c_\varepsilon=\sup\left|g_{z,v}[\lambda_1,\ldots,\lambda_k]\right|<\infty$
where the supremum is taken over $z\in D_{-\frac\varepsilon2}$, all unit vector $v$ in $\cc^n$, all positive integer $k$ such that $k\leq \#\Lambda_{z,v}$, all $\lambda_1,\ldots,\lambda_k\in\Lambda_{z,v}$ pairwise distinct.\\
Then $g$ as a smooth extension on $D_{-\varepsilon}$ bounded by $c_\varepsilon$ up to a multiplicative constant uniform with respect to $g$.
\end{lemma}
\pr We proceed locally  and glue all the extension. Since the only problems occur when we are near a singularity we consider $z_0$ a singularity of $X$ and we choose an orthonormal basis $w_1,w_2$ such that $\pi_0$,  the orthogonal projection on the hyperplane orthogonal to $w_2$ passing through $z_0$, is a $k_0$ sheeted covering of $X$ in a neighborhood ${\cal U}_0\subset D$ of $z_0$.\\
For $z_1\neq 0$, we denote by $\lambda_1(z_1),\ldots,\lambda_{k_0}(z_1)$ the pairwise distinct complex number such that for $k=1,\ldots, k_0$, $z_0+z_1w_1+\lambda_k(z_1) w_2$ belongs to $X$. We set for $z=z_0+z_1w_1+z_2w_2$, $z_1\neq 0$ :
$$\tilde g_0(z)=\tilde g_0(z_0+z_1w_1+z_2w_2)=\sum_{k=1}^{k_0}\prod_{\over{l=1}{l\neq k}}^{k_0} \frac {z_2-\lambda_l(z_1)}{\lambda_k(z_1)-\lambda_l(z_1)} g(z_0+z_1w_1+\lambda_k(z_1) w_2).$$
By construction, $\tilde g_0(z)=g(z)$ for all $z\in X\cap {\cal U}_0$, $z\neq z_0$. We denote by $\Delta_0$ the complex line passing through $z_0$ and supported by $w_2$.\\
Since $z_0$ is an isolated singularity of $X$, away from $0$, the $\lambda_j$ depend locally holomorphicaly from $z_1$ and thus $\tilde g_0$ is holomorphic on ${\cal U}_0\setminus \Delta_0$.\\
Since the divided differences are bounded on $D_{-\frac\varepsilon2}$ by $c_\varepsilon $, $\tilde g_0$ is bounded on $ {\cal U}_0\setminus \Delta_0$ by $c_\varepsilon$ up to a uniform multiplicative constant and thus $\tilde g_0$ is holomorphic and bounded on ${\cal U}_0$.\qed

Combining Theorems \ref{th0}, \ref{th5}, Lemma \ref{lemma1} and corollary \ref{th2} we get the two following theorems.
\begin{theorem}\label{th8}
 For $X$ and $D$ in $\cc^2$, let $g$ be a weakly holomorphic function in $\cc^2$ such that $c_\infty(g)<\infty$. Then $g$ has a holomorphic extension $G$ which belong to $BMO(D)$ such that $\|G\|_{BMO(D)}\leqs c_\infty(g)$.
\end{theorem}
\begin{theorem}\label{th9}
For $X$ and $D$ in $\cc^2$, let $g$ be a weakly holomorphic function in $\cc^2$ such that $c^{(q)}_{\kappa,(z_j)_{j\in\nn}}(g)<\infty$ and $c_\varepsilon<\infty$. Then $g$ has a holomorphic extension $G$ which belongs to $L^q(D)$ such that $\|G\|_{L^q(D)}\leqs c^{(q)}_{\kappa,(z_j)_{j\in\nn}}(g)+c_\varepsilon(g)$.
\end{theorem}

\section{Examples}\mlabel{section6}
\begin{example}[$BMO$ extension]
Let $D$ be the ball of radius 1 and center $(1,0)$ in $\cc^2$. We choose $\rho(z)=|z_1-1|^2+|z_2|^2-1$ as a defining function for $D$. For $\alpha_1,\alpha_2,\ldots,\alpha_k\in\cc$ pairwise distinct we set $v_i=(-\overline{\alpha_i},1)$. We  denote by  $P_i$ the plane orthogonal to $v_i$ passing through the origin and we set $\Delta_i=P_i\cap D$ and $X=\cup_{i=1}^k P_i$. Let also $g_1,\ldots, g_k$ be $k$ bounded holomorphic functions on $\Delta$, the unit disc in $\cc$. Since $\Delta_i=\{(z_1,z_2)\in\cc^2,\ z_2=\alpha_iz_1\text{ and } |z_1-(1+|\alpha_i|^2)^{-1}|<(1+|\alpha_i|^2)^{-1}\}$, the function
$$\app g {X\cap D} \cc {(z_1,z_2)}{g_i(z_1(1+|\alpha _i|^2)-1)}$$
is well defined, bounded and holomorphic on $X\cap D$. Question~: Under which conditions does $g$ have a $BMO$ holomorphic extension on the domain $D$ ?
\end{example}
In order to answer this question, we will try to find an upper bound for $c^{(\infty)}_{\kappa,\varepsilon_0}(g)$. Let $\zeta=(\zeta_1,\zeta_2)$ be a point in $D\setminus D_{-\varepsilon_0}$, let $z^*_1\in\cc$ be such that $|z^*_1|<\kappa|\rho(\zeta)|$ and let $\lambda _1,\ldots, \lambda _l$ be complex numbers pairwise distinct belonging to $\Lambda_{\alpha+z^*_1\eta_\zeta,v_\zeta}$. Perhaps after renumbering, we assume that $\zeta+z^*_1\eta_\zeta+\lambda_iv\zeta$ belongs to $\Delta_i$ for all $i$. Moreover, if $\zeta$ is sufficiently near the origin, we can also assume that $v_\zeta$ does not belong to any of the plane $P_i$.\\
We have
\begin{eqnarray*}
g_{\zeta+z^*_1\eta_\zeta,v_\zeta}[\lambda _1,\ldots,\lambda _l]&=&\sum_{i=1}^l\frac1{\prod_{\over{j=1}{j\neq i}}^l(\lambda_i-\lambda_j)} g_i\left((\zeta_1+z^*_1\eta_{\zeta,1}+\lambda_iv_{\zeta,1}) (1+|\alpha _i|^2)-1\right) .
\end{eqnarray*}
For $m=i,j$, $\lambda_m$ satisfies the following equalities
$$\zeta_2+z^*_1\eta_{\zeta,2}+\lambda_mv_{\zeta,2}=\alpha_m(\zeta_1+z^*_1\eta_{\zeta,1}+\lambda_lv_{\zeta,1}),\qquad m=i,j$$
which yield $\lambda_i-\lambda_j=(\alpha_i-\alpha_j)(\zeta_1z^*_1\eta_{\zeta,1}+\lambda_iv_{\zeta,1})+\alpha_j(\lambda_i-\lambda_j)v_{\zeta,1}$
and so 
$$|\lambda_i-\lambda_j|\cdot |v_{\zeta,2}-\alpha_jv_{\zeta,1}|=|\alpha_i-\alpha_j|\cdot|\zeta_1+z^*_1\eta_{\zeta,1}+\lambda_iv_{\zeta,1}|.$$
We show that $|\zeta_1+z^*_1\eta_{\zeta,1}+\lambda_iv_{\zeta,1}|\geqs |\zeta_1|$.\\
First, we have $|z_1^*|\leq \kappa|\rho(\zeta)|$ and since $\zeta$ belongs to $D$, $|\rho(\zeta)|\leqs|\zeta_1|$ so $|z^*_1|\leqs \kappa|\zeta_1|$.\\
Secondly, $|v_{\zeta,1}|\eqs \left|\diffp{\rho}{\zeta_2}(\zeta)\right|\eqs |\zeta_2|$ and since $\zeta$ belongs to $D$, $|\zeta_2|\leqs |\zeta_1|^{\frac12}$. Since $|\lambda_i|\leq 3\kappa|\rho(\zeta)|^{\frac12}\leq |\zeta_1|^{\frac12}$, we get $|\lambda_i v_{\zeta,1}|\leqs \kappa|\zeta_1|$.\\
Thus provided $\kappa$ is small enough, $|\lambda_i-\lambda_j|\geqs |\zeta_1|$ and
\begin{eqnarray*}
\left|g_{\zeta+z^*_1\eta_\zeta,v_\zeta}[\lambda _1,\ldots,\lambda _l]\right|&\leqs&\frac1{|\zeta_1|^{l-1}} \sum_{i=1}^l \left| g_i\left((\zeta_1+z^*_1\eta_{\zeta,1}+\lambda_iv_{\zeta,1}) (1+|\alpha _i|^2)-1\right)\right|.
\end{eqnarray*}
Since $\tau(\zeta,v_\zeta,|\rho(\zeta)|)\leqs |\zeta_1|^{\frac12}$, if we assume that there exists $c\in\cc$ and $C>0$ such that for all $i$, $|g_i(z+1)-c|\leq C |z|^{\frac{l-1}2}$ for all $z$ near the origin of $\cc$, we get
$$\tau(\zeta,v_\zeta,|\rho(\zeta)|)^{l-1}\left|g_{\zeta+z^*_1\eta_\zeta,v_\zeta}[\lambda _1,\ldots,\lambda _l]\right|\leqs C.$$
So $c^{(\infty)}_{\kappa,\varepsilon_0}(g)$ is finite and Lemma \ref{lemma2} and Theorem \ref{th0} implies that $g$ admits a $BMO$-holomorphic extension on $D$.
\par\medskip
This is in general the best result we can get. For example, let $\alpha$ be a real number and let $g_i$ be the function defined on the unit disc of $\cc$ by $g_i(z)=(1+z)^\alpha$, $i=1,\ldots, k$. Let $x$ be a small positive real number and let $\zeta$ in $D$  be the point $(x,0)$. We have $\eta_\zeta=(1,0)$, $v_\zeta=(0,1)$, $\tau(\zeta,v_\zeta,|\rho(\zeta)|)\eqs x^{\frac12}$, $(x,\alpha_ix)$ belongs to $\Delta_i$ if $x$ is sufficiently small and
$$g_{\zeta,v_\zeta}[\alpha_1x,\ldots, \alpha_kx]=\sum_{i=1}^k \frac 1{x^{k-1}\prod_{\over{j=1}{j\neq i}} (\alpha_i-\alpha_j) }\left(x(1+|\alpha_i|^2)\right)^\alpha.$$
Therefore if $\alpha<\frac {k-1}2$, $\tau(\zeta,v_\zeta,|\rho(\zeta)|)^{k-1} |g_{\zeta,v_\zeta}[\alpha_1x,\ldots, \alpha_kx]|$ is unbounded when $x$ goes to $0$. So $c_\infty(g)$ is not finite and Theorem \ref{th1} implies that $g$ does not admit a holomorphic extension bounded on $D$.

\begin{example}[$L^2$-extension in $\cc^2$]
 Again let $D$ be the ball of radius 1 and center $(1,0)$ in $\cc^2$ and for any positive odd integer $q$, let $X$ be the analytic set $X=\{z\in\cc^2,\ z^q_1=z^2_2\}$. Then all $g$ holomorphic and bounded on $X\cap D$ has a $L^2$ holomorphic extension on $D$ if and only if $q=1$ or $q=3$.
\end{example}
When $q=1$, $X$ a manifold and there is nothing to do.

When $q=3$, $X$ has a singularity at the origin. We will prove that the assumptions of Theorem \ref{th5} are satisfied for any $\kappa$-covering provided $\kappa$ is small enough. To check these hypothesis,
we set $\rho(z)=|z_1-1|^2+|z_2|^2-1$, we fix a holomorphic square root $\alpha$ in $\cc\setminus(-\infty,0]$ and we prove the following facts. The first one gives a relation between the distance from $z\in X\cap D$ to $z+\lambda v\in X\cap D$ and the coordinates of $z$.
\begin{fact}\label{facta}
 Let $\kappa $ be a sufficiently small positive real number, let $K$ be a large positive real number, let $z=(z_1,z_2)$ be a point in $D\cap X$, let $v=(v_1,v_2)$ be a unit vector of $\cc^2$ such that $|v_1|\leq K|z_1|^{\frac12}$ and let $\lambda$ be a complex number such that $z+\lambda v$ belongs to $X\cap D$ and $|\lambda|\leq 4\kappa |\tau(z,v,|\rho(z)|)$.\\
Then, if $\kappa $ is small enough, we have $|\lambda|\geqs |z_1|^{\frac q2}$, $|z_1|\leqs |\rho(z)|^\frac 1q$ and $|z_2|\leqs |\rho(z)|^{\frac12}$ each time uniformly with respect to $z$, $\kappa$ and $v$.
\end{fact}
\begin{remark}
 The assumption  $|v_1|\leq K|z_1|^{\frac12}$ means that $v$ is ``nearly'' tangential to $bD_{\rho(z)}$.
\end{remark}

\pr We first prove that $|\lambda|\geqs |\rho(z)|^{\frac q2}$.  Since $v$ is transverse to $X$, without restriction we assume that $z=(z_1,\alpha(z_1)^q)$ and that $z+\lambda v=(z_1,-\alpha(z_1+\lambda v_1)^q)$. Therefore we have 
$$|\lambda|\geq|\alpha^q(z_1)+\alpha^q(z_1+\lambda v_1)|\\
\geq2|z_1|^{\frac q2}-|\alpha^q(z_1)-\alpha^q(z_1+\lambda v_1)|.$$
The mean value theorem gives
$$|\alpha^q(z_1)-\alpha^q(z_1+\lambda v_1)|\leqs |\lambda||v_1| \sup_{\zeta\in[z_1,z_1+\lambda v_1]} \left|\diffp {\alpha^q}{\zeta}(\zeta)\right|.$$
For all $\zeta\in[z_1,z_1+\lambda v_1]$, we have $|\zeta|\leqs |z_1|$, and so, provided $\kappa $ is small enough, we get $|\lambda|\geq |z_1|^{\frac q2}$.
Now, since $|\lambda|\leq 4\kappa |\rho(z)|^{\frac12}$, we get $|z_1|\leqs |\rho(z)|^{\frac1q}$ and $|z_2|\leqs|\rho(z)|^{\frac12}$.\qed

As previously, we denote by $\eta_\zeta$ the outer unit normal to $bD_{\rho(\zeta)}$ at $\zeta$ and by $v_\zeta$  a tangent vector to $bD_{\rho(\zeta)}$ at $\zeta$. The second fact gives some kind of  uniformity of Fact \ref{facta} on a Koranyi ball.
\begin{fact}\label{factb}
 Let $\kappa $ be a sufficiently small positive real number, let $\zeta$ be a point in $D$, let $z=\zeta+z_1^*\eta_\zeta+z_2^*v_\zeta$ be a point in ${\cal P}_{4\kappa |\rho(\zeta)|}(\zeta)\cap D\cap X$ and let $\lambda$ be a complex number such that $z+\lambda v_\zeta$ belongs to $X\cap D\cap{\cal P}_{4\kappa |\rho(\zeta)|}(\zeta)$.\\
Then $|\lambda|\geqs |\zeta_1|^{\frac q2}$, $|\zeta_2|\leqs|\rho(\zeta)|^{\frac12}$ and $|\zeta_1|\leqs|\rho(\zeta)|^{\frac1q}$ uniformly with respect to $z$, $\zeta$ and $\lambda $.
\end{fact}
\pr
We want to apply Fact \ref{facta}, so we first have to check that $|v_{\zeta,1}|\leqs |z_1|^{\frac12}$, uniformly with respect to $z$ and $\zeta$.\\
On the one hand we have $|v_{\zeta,1}|\eqs\left|\diffp{\rho}{\zeta_2}(\zeta)\right|\eqs|\zeta_2|\leqs |\zeta_1|^{\frac 12}$.\\
On the other hand $z_1=\zeta_1+z^*_1\eta_{\zeta,1}+z^*_2v_{\zeta,1}$ thus 
\begin{eqnarray*}
 |\zeta_1|&\leq& |z_1^*|+  |z_2^*||v_{\zeta,1}|+|z_1|\\
&\leqs&  \kappa |\rho(z)|+ \kappa |v_{\zeta,1}|^2+|z_1|\\
&\leqs& |z_1| + \kappa |v_{\zeta,1}|^2.
\end{eqnarray*}
Therefore, if $\kappa $ is small enough,  $|v_{\zeta,1}|\leqs |z_1|^{\frac12}$ and $|\zeta_1|\leqs|z_1|$. Therefore we can apply Fact \ref{facta} which gives 
$|\lambda|\geqs |z_1|^{\frac q2}$ and since  $|z_1|\geqs |\zeta_1|$ the first inequality is proved. The third inequality follows from the first one and from the fact that $|\lambda|\leqs|\rho(\zeta)|^{\frac12}$.\\
Fact \ref{facta} also gives $|z_2|\leqs |\rho(z)|^{\frac12}$ and since $|\rho (\zeta)|\eqs|\rho (z)|$, we have
$$|\zeta_2|\leqs |\zeta_2-z_2|+|z_2|\leqs |\rho(\zeta)|^{\frac12}+|\rho(z)|^{\frac12}\leqs|\rho(\zeta)|^{\frac12}.\qed$$
Now we check the assumptions of Theorem \ref{th5} and for any $\kappa$-covering, $\kappa>0$ sufficiently small, and any function $g$ bounded on $X\cap D$ we  prove that $c^{(2)}_{\kappa ,(\zeta_j)_{j\in\nn}}(g)\leqs\|g\|_{L^\infty(D\cap X)}$, uniformly with respect to $g$.\\
Let ${\cal U}_0$ be a neighborhood of the origin, let $c$, $\varepsilon_0$ and $\kappa$ be small positive real numbers and let ${\cal P}_{\kappa|\rho(\zeta^{(k)}_j)|}(\zeta^{(k)}_j)$, $k\in\nn$, $j\in\{1,\ldots, n_k\}$ be a $\kappa$-covering of $D\cap {\cal U}_0$ such that for all $k$ and all $j$, the point $\zeta_j^{(k)}$ belongs to $bD_{-(1-c\kappa)^k\varepsilon_0}$. We assume that $\kappa$ is so small that Fact \ref{factb} holds true and we set $\tilde \kappa=1-c\kappa$.

For all $\zeta\in D$, the following inequality holds
\begin{eqnarray*}
 |\rho(\zeta)|\int_{|z^*_1|<4\kappa |\rho (\zeta)|}\sum_{\lambda\in\Lambda_{\zeta +z^*_1\eta_\zeta ,v_\zeta }} \left|g_{\zeta +z^*_1\eta_\zeta ,v_\zeta }[\lambda]\right|^2dV(z_1)&\leqs &\|g\|_{L^\infty(X\cap D)}^2 |\rho (\zeta)|^3.
\end{eqnarray*}
This means that the corresponding estimate for $\zeta^{(k)}_j$ does not depend on $j$ and since we will add these bound for all $k$ and $j=1,\ldots, n_k,$ we will also need an upper bound for $n_k$. For any non negative integer $k$, we denote by $\sigma_k$ the area measure on $bD_{-\tilde \kappa^k\varepsilon_0}$. Since ${\cal P}_{\kappa|\rho(\zeta^{(k)}_j)|}(\zeta^{(k)}_j)$ is a $\kappa$-covering, for all $k$ we have as in the proof of Proposition \ref{propmax}
\begin{eqnarray*}
 \sigma_k\left (bD_{\tilde \kappa ^k\varepsilon _0}\right)&\geq&\sigma_k\left(bD_{\tilde \kappa ^k\varepsilon _0}\cap \cup_{j=1}^{n_k}{\cal P}_{\kappa |\rho (\zeta_j^{(k)})|}((\zeta_j^{(k)}))\right)\\
 &\geq& \sum_{j=1}^{n_k} \sigma_k\left(bD_{\tilde \kappa ^k\varepsilon _0}\cap {\cal P}_{\frac c{c_1}\kappa |\rho (\zeta_j^{(k)})|}((\zeta_j^{(k)}))\right)\\
 &\geqs& n_k \left(\tilde \kappa^k\varepsilon_0\right)^{2}.
\end{eqnarray*}
Therefore $n_k\leqs (\tilde \kappa^k\varepsilon_0)^{-2}$ and  we have uniformly with respect to $g$
\begin{eqnarray*}
\lefteqn{\sum_{k=0}^\infty \sum_{j=1}^{n_k} |\rho(\zeta_j^{(k)})|\int_{|z_1^*|<4\kappa |\rho (\zeta_j^{(k)})|}\sum_{\lambda\in\Lambda_{\zeta +z^*_1\eta_{\zeta^{(k)}_j} ,v_{\zeta^{(k)}_j} }}\left|g_{\zeta_j^{(k)}+z^*_1\eta_{\zeta^{(k)}_j},v_{\zeta^{(k)}_j}}[\lambda]\right|^2dV(z^*_1)}\\
&&\hskip230pt\leqs\|g\|_{L^\infty(X\cap D)}^2
\sum_{k=0}^{\infty} n_k \left(\tilde\kappa^k\varepsilon _0\right)^3\\
&&\hskip230pt\leqs\|g\|_{L^\infty(X\cap D)}^2.
\end{eqnarray*}
Now we handle the case of divided differences of order 2. We set 
$$I(\zeta)=|\rho(\zeta)|^2\int_{|z^*_1|<4\kappa |\rho (\zeta)|} \sum_{
\over{\lambda_1,\lambda _2\in \Lambda _{\zeta +z^*_1\eta_\zeta,v_\zeta}}{\lambda _1\neq\lambda _2}}
\left|g_{\zeta +z^*_1\eta_\zeta,v_\zeta}[\lambda_1,\lambda_2] \right|^2dV(z^*_1)$$
and we aim to prove that $\sum_{k=0}^{+\infty}\sum_{j=1}^{n_k} I(\zeta_j^{(k)})\leqs\|g\|_{L^\infty(X\cap D)}$.\\
If for all complex number $z^*_1$ such that $|z^*_1|\leq \kappa |\rho (\zeta)|$ we have $\#\Lambda _{\zeta +z^*_1\eta_\zeta,v_\zeta}<2$, then $I(\zeta)=0$. Otherwise
Fact \ref{factb} implies that  $|\zeta_2|\leq K (\tilde \kappa \varepsilon _0)^{\frac12}$ for some $K>0$ and that 
$|\lambda _1-\lambda _2|\geqs|\zeta_1|^{\frac32}$
for all $\lambda_1,\lambda_2$ distinct in $\Lambda _{\zeta+z^*_1\eta_\zeta,v_\zeta}$, $z_1^*\in\cc$ such that $|z_1^*|\leq \kappa|\rho (\zeta)|$. Therefore, for all such $\zeta$, we have
\begin{eqnarray}
 I(\zeta) \leqs|\rho(\zeta )|^2\int_{|z_1^*|<4\kappa |\rho (\zeta)|} \frac{\|g\|_{L^\infty(D\cap X)}} {|\zeta_1|^3}dV(z^*_1)\leqs \|g\|_{L^{\infty}(X\cap D)} \frac{|\rho(\zeta)|^4}{|\zeta_1|^3}\label{eq31}
\end{eqnarray}
Thus, when we denote by $Z^{(k)}$ the set
$$Z^{(k)}=\{j\in\nn,\ \exists z^*_1\in\cc,\ |z^*_1|<\kappa|\rho (\zeta_j^{(k)})| \text{ and } \#\Lambda _{\zeta_j^{(k)} +z^*_1\eta_{\zeta_j^{(k)}},v_{\zeta_j^{(k)}}}=2\},$$
we have to estimate the sum $\sum_{k=0}^{+\infty} \sum_{j\in Z^{(k)} } \frac{(\tilde \kappa^k \varepsilon _0)^4}{|\zeta^{(k)}_{j,1}|^3}$.\\
We write $Z^{(k)}$ as  $Z^{(k)}=\cup_{i=1}^{\infty} Z^{(k)}_i$ where $Z^{(i)}_{k}=\{j\in Z^{(k)},\ i\tilde\kappa^k \varepsilon _0\leq |\zeta^{(k)}_{j,1}|< (i+1)\tilde\kappa ^k\varepsilon _0\text{ and } |\zeta^{(k)}_{j,2}|\leq K (\tilde \kappa \varepsilon _0)^{\frac12} \}$
and we look for an upper bound of $\#Z_{i}^{(k)}$.
We have
$$\sigma _k(bD_{-\tilde\kappa ^k\varepsilon _0}\cap \{z,\ \frac12i\tilde \kappa ^k\varepsilon _0\leq|z_1|\leq 2(i+1)\tilde \kappa ^k\varepsilon _0\text{ and } |z_2|\leq 2K(\tilde \kappa ^k\varepsilon _0)^{\frac12}\})\eqs (\tilde \kappa ^k\varepsilon _0)^2$$
and, if $\kappa$ is small enough :
\begin{eqnarray*}
 \lefteqn{\sigma _k(bD_{-\tilde\kappa ^k\varepsilon _0}\cap \{z,\ \frac12i\tilde \kappa ^k\varepsilon _0\leq|z_1|\leq 2(i+1)\tilde \kappa ^k\varepsilon _0\text{ and }  |z_2|\leq K(\tilde \kappa ^k\varepsilon _0)^{\frac12}\})}\\
 &&\hskip200pt\geqs \sigma_k(\cup_{j\in Z_{i}^{(k)}} {\cal P}_{\kappa|\rho (\zeta^{(k)}_j)|}(\zeta_j^{(k)}) \cap bD_{-\tilde\kappa ^k\varepsilon _0})\\
 &&\hskip200pt\geqs \#Z_i^{(k)} \cdot (\tilde \kappa ^k\varepsilon _0)^2.
\end{eqnarray*}
These last two inequalities imply that
$\# Z_i^{(k)}$ is bounded by a constant which does not depend from $i$ nor from $k$. \\
For $j\in Z_0^{(k)}$, since $|\zeta^{(k)}_{j,1}|\geqs |\rho (\zeta_j^{(k)})|$, Inequality (\ref{eq31}) yields $I(\zeta_j^{(k)})\leqs \tilde\kappa ^k\varepsilon _0\|g\|_{L^\infty(X\cap D)}$ thus $$\sum_{k=0}^{+\infty} \sum_{j\in Z_0^{(k)}} I(\zeta^{(k)}_j)\leqs \|g\|_{L^{\infty}(X\cap D)}.$$
For $i>0$, we use directly (\ref{eq31}) which gives 
$$
 \sum_{i=1}^{+\infty} \sum_{k=0}^{+\infty}\sum_{j\in Z_i^{(k)}}I(\zeta^{(k)}_j) 
\leqs \|g\|_{L^{\infty}(X\cap D)}
 \sum_{k=0}^{+\infty} \sum_{i=1}^{+\infty} \frac{(\tilde\kappa ^k\varepsilon _0)^4}{(i\tilde \kappa ^k\varepsilon _0)^3}\\
 \leqs\|g\|_{L^{\infty}(X\cap D)}.
$$
This ends to prove that $c^{(2)}_{\kappa,(\zeta_j^{(k)})_{k\in\nn,j\in\{1,\ldots,n_k\}}}$ is finite and Theorem \ref{th5} now implies that $g$ admits  a $L^2$-holomorphic extension on $D$.
\par\medskip
Now, for $q\geq 5$, we consider $g$ defined for $z$ in $X$ by $g(z)=\frac{z_2}{z_1^{\frac q2}}$. The function $g$ is holomorphic and bounded on $X$  because $|z_2|=|z_1|^{\frac q2}$ for all $(z_1,z_2)\in X$ but we will see that $g$ does not admits a $L^2$-holomorphic extension on $D$.\\
For $\varepsilon _0, \kappa,c>0$ small enough we set $\tilde \kappa =1-c\kappa $ and we denote by $\zeta^{(k)}_0=(x_k,0)$ the point of $\cc^2$ such $\rho(\zeta_0^{(k)})= -\tilde \kappa ^k\varepsilon _0$. We have $x_k\eqs \tilde \kappa^k\varepsilon _0$ uniformly with respect to $k$, $\kappa$ and $\varepsilon _0$. We complete the sequence $(\zeta_0^{(k)})_{k\in\nn}$ so as to get a $\kappa$-covering ${\cal P}_{\kappa|\rho(\zeta^{(k)}_j)|}(\zeta^{(k)}_j)$, $k\in\nn$ and $j\in\{0,\ldots,n_k\}$, of a neighborhood of the origin. We set $w_1=(1,0)$ and $w_2=(0,1)$. For all $k$, $\eta_{\zeta^{(k)}_0}=w_1$, $v_{\zeta^{(k)}_0}=w_2$ and, for all $z_1$, we have $\Lambda_{\zeta^{(k)}_0+z_1w_1,w_2}=\{(z_1+\tilde\kappa^k\varepsilon_0)^{\frac q2},-(z_1+\tilde\kappa^k\varepsilon_0)^{\frac q2}\}$. So, if $\kappa$ is small enough,  for all $k$ we have
\begin{eqnarray*}
 \lefteqn{|\rho(\zeta _0^{(k)})|^2 \int_{|z_1|<4\kappa |\rho (\zeta _0^{(k)})|} \left|g_{\zeta_0^{(k)}+z_1w_1,w_2}\left[(z_1+\tilde\kappa^k\varepsilon_0)^{\frac q2},-(z_1+\tilde\kappa^k\varepsilon_0)^{\frac q2}\right]\right|^2dV(z_1)}\\
&& \hskip 160pt \geqs (\tilde \kappa ^k\varepsilon _0)^2\int_{|z_1|<4\kappa |\rho (\zeta_0^{(k)})|} \frac1{|z_1+\tilde\kappa^k\varepsilon_0|^q}dV(z_1)\\
&&\hskip 160pt \geqs (\tilde \kappa ^k\varepsilon _0)^{4-q}.
\end{eqnarray*} 
Since for $q\geq 5$ the series $\sum_{k\geq 0} (\tilde \kappa ^k\varepsilon _0)^{4-q}$ diverges $c^{(2)}_{\kappa,(\zeta^{(k)}_j)_{k\in\nn,j\in\{0,\ldots, n_k\}}}(g)$ is not finite and so  Theorem \ref{th4} implies that $g$ does not have a $L^2$ holomorphic extension on $D$.
\begin{example}[The example of Diederich-Mazzilli]
Let $B_3$ be the unit ball of $\cc^3$, $X=\{z=(z_1,z_2,z_3)\in\cc^3\ : \ z_1^2+z^q_2=0\}$ where $q\geq 10$ is an uneven integer and define the holomorphic function $f$ on $\cc^3$ by 
$$f(z)=\frac{z_1}{(1-z_3)^{\frac q4}}.$$
Then $f$ is bounded on $X\cap B_3$ and has no $L^2$ holomorphic extension on $B_3$.
\end{example}
This was shown in \cite{DiMa0} by Diederich and the second author. We will prove this result here with Theorem \ref{th4}.

We set $\rho(\zeta )=|\zeta _1|^2+|\zeta_2|^2+|\zeta_3|^2-1$, and we denote by $w_1,w_2, w_3$ the canonical basis of $\cc^3$.
For all non negative integer $j$ and $\varepsilon _0,c$ and $\kappa $ small suitable constants for $X$ and $B_3$, we define $\tilde\kappa =(1-c\kappa )$. For any integer $j$, we denote by $\zeta_j=(0,0,\zeta_{j,3})$ the point of $\cc^3$ such that $\zeta_{j,3}$ is real and satisfies $\rho(\zeta_j)=-\tilde \kappa^j\varepsilon _0$. The point $\zeta _j$ can be chosen at the first step of the construction of a $\kappa $-covering of $X\cap D$ in a neighborhood of $(0,0,1)$ and so the Koranyi balls 
${\cal P}_{\kappa|\rho(\zeta_j)|}(\zeta_j)$, $j\in\nn$, are extract from a $\kappa$-covering.
For all $j$ we have
\begin{eqnarray*}
 |\rho(\zeta _j)|^2\int_{\over{|z_2|<(4\kappa|\rho(\zeta _j)|)^{\frac12}}{|z_3-\zeta _{j,3}|<4\kappa|\rho (\zeta _j)|} }\left|f_{\zeta_j +z_2 w_2+z_3w_3,w_1} \left[z_2^{\frac q2},-z_2^{\frac q2}\right]\right|^2dV(z_2,z_3)&\geqs & \tilde\kappa ^{j(5-\frac q2)}
\end{eqnarray*}
and thus when $q\geq 5$, 
\begin{eqnarray*}
\sum_{j=0}^{+\infty}|\rho(\zeta _j)|^2\int_{\over{|z_2|<(4\kappa|\rho(\zeta _j)|)^{\frac12}}{|z_3-\zeta _{j,3}|<4\kappa|\rho (\zeta _j)|} }\left|f_{\zeta_j +z_2w_2+z_3w_3,w_1} \left[z_2^{\frac q2},-z_2^{\frac q2}\right]\right|^2dV(z_1,z_3)&=&+\infty.
\end{eqnarray*}
Theorem \ref{th4} then implies that $f$ does not have a $L^2$ holomorphic extension on $B_3$.

\end{document}